\documentclass[a4paper,12pt,english]{amsart}
\usepackage[a4paper]{geometry}
\usepackage[utf8]{inputenc}
\usepackage[T1,OT1]{fontenc}
\newcommand{\ogonek}[1]{{\fontencoding{T1}\selectfont\k{#1}}}
\usepackage{cite}
\usepackage{url}
\usepackage{xcolor}
\usepackage{graphicx}
\usepackage{soul}
\usepackage[normalem]{ulem}
\usepackage{mathtools}  
\mathtoolsset{showonlyrefs}
\usepackage{hyperref}

\usepackage[foot]{amsaddr}
\usepackage{enumerate}

\usepackage{amssymb}
\usepackage{array}
\allowdisplaybreaks[4]
\usepackage{bbm}

\newtheorem{theorem}{Theorem}[section]
\newtheorem{lemma}[theorem]{Lemma}
\newtheorem{definition}[theorem]{Definition}

\newtheorem{proposition}[theorem]{Proposition}
\newtheorem{remark}[theorem]{Remark}

\numberwithin{equation}{section}

\newcommand{\Bcal} {{\mathcal B}}

\newcommand{\Ecal} {{\mathcal E}}
\newcommand{\Fcal} {{\mathcal F}}

\newcommand{\Mcal} {{\mathcal M}}

\newcommand{\Wcal} {{\mathcal W}}

\newcommand{\R}{\mathbb{R}}

\renewcommand{\P}{\mathbb{P}}

\newcommand{\EE}{\mathbb{E}}

\renewcommand{\epsilon}{\varepsilon}
\newcommand{\eps}{\varepsilon}

\newcommand{\Lip}{\operatorname{Lip}}

\newcommand{\ud}{\mathrm{d}}
\newcommand{\e}{\varepsilon}

\usepackage{graphicx}
\usepackage{tikz}
\usepackage{atbegshi}

\AtBeginShipoutFirst{%
    \begin{tikzpicture}[remember picture, overlay]
        \node[anchor=north west, xshift=1in, yshift=-1.5cm] at (current page.north west) {%
            \includegraphics[width=2cm]{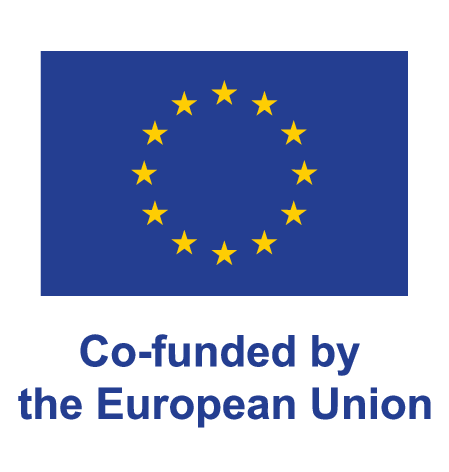}%
        };
    \end{tikzpicture}%
}

\title[Mixing times for the stochastic $p$-Laplace equation]{Mixing times for the stochastic $p$-Laplace equation}
\author[G. Barrera]{Gerardo Barrera}
\email{gerardo.barrera.vargas@tecnico.ulisboa.pt}

\address[GB]{Center for Mathematical Analysis\\ Geometry and Dynamical Systems\\ Instituto
Superior T\'ecnico\\ Universidade de Lisboa\\ Av. Rovisco Pais\\ 1049-001 Lisboa\\
Portugal}

\author[J. M. T\"olle]{Jonas M. T\"olle}
\email{jonas.tolle@aalto.fi}

\address[JMT]{Aalto University\\
Department of Mathematics and Systems Analysis\\
PO Box 11100 (Otakaari 1, Espoo)\\
00076 Aalto\\
Finland}

\date{\today}

\begin{document}
\begin{abstract}
We review quantitative results on the convergence to equilibrium of the stochastic $p$-Laplace equation, $p>1$, on a bounded domain in $\R^d$ with additive, possibly degenerate Wiener noise. From these results we derive upper bounds on $\eps$-mixing times in different distances, depending on the initial data. We show that the Wasserstein distance to equilibrium is non-increasing in time and that, for $p>\frac{2d}{d+2}$, the invariant measure has moments of all orders. For $1<p<2$ with $p>\frac{2d}{d+2}$ and possibly degenerate noise, the $\eps$-mixing times in every Wasserstein distance $\mathcal W_r$, $r\ge1$, are $O(\eps^{-\delta})$ as $\eps\downarrow0$ for every $\delta>0$. Two tables collect the convergence rates and the mixing-time bounds.
\end{abstract}
\maketitle
{\footnotesize\tableofcontents}

\section{Introduction}

It is well known that no information about the rate of convergence can be deduced from Birkhoff's ergodic theorem alone, see~\cite{K:78}. For quasi-periodic and almost periodic systems and sufficiently regular observables, however, suitably weighted Birkhoff averages converge exponentially fast or faster than any power of time, see~\cite{TL:24,TL:26}. Von Neumann's mean ergodic theorem~\cite{vN:32} likewise yields no rate of convergence in general. Rates of convergence in ergodic theorems are discussed in~\cite{K:96}.
In many cases, physically relevant quantitative results for ergodic systems can be obtained by precise analysis of the underlying dynamical system.

Ergodic properties of stochastic partial differential equations (SPDEs) with applications to Kolmogorov operators on infinite-dimensional spaces have been systematically studied in~\cite{DPZ:96,BDP:06,DPGZ:92,DP:04}. The literature on this topic has developed extensively.
In the recent decade, the approach by coupling has turned out to be most fruitful, see~\cite{KS:18,BKS:20,W:15-1,K:18,Eberle}.

Mixing times have been studied in~\cite{DIA:96,LPW:09} in the context of Markov chains.
They quantify the time needed, from a given initial datum, to come $\eps$-close to equilibrium in a chosen distance. Mixing times for stochastic differential equations (SDEs) have been studied in~\cite{BHP:24,BP:20,BCL:22}, among others. Mixing times for SPDEs have been studied by the first author with H\"{o}gele and Pardo in~\cite{BHP:23} and by the present authors in~\cite{BT:24}.

This work aims to give an overview of the ergodic rates of the stochastic $p$-Laplace equation, see Subsection~\ref{subsec:pLap} below. The $p$-Laplace equation is an important example of a nonlinear parabolic diffusion equation with $p$-growth that can be treated by variational methods. It serves as a toy model for power-law fluids~\cite{B:15} and is discussed in detail in the monograph~\cite{DB:93}. Our contributions are described in Subsection~\ref{subsec:contrib}.

\subsection{Motivating example}\label{subsec:motivation}

Let $\{\beta_t\}_{t\ge 0}$ be a standard Brownian motion on $\R$ defined on a filtered probability space $(\Omega,\Fcal,\{\Fcal_t\}_{t\ge 0},\P)$ satisfying the usual conditions.
Consider the unique strong solution $\{X_t^x\}_{t\ge 0}$ of the following It\^o stochastic differential equation (SDE) with $(p-1)$-homogeneity in the drift, $p>1$,
\begin{equation}\label{eq:SDE}
\ud X_t^x=-(X_t^x)^{[p-1]}\,\ud t+\ud \beta_t,\quad X_0^x=x\in\R,
\end{equation}
where we use the notation $z^{[\alpha]}:=z|z|^{\alpha-1}$ for $z\in\R\setminus\{0\}$ and $\alpha>0$, and we set $0^{[\alpha]}:=0$.
We write $(P_t)_{t\ge0}$ for the transition semigroup of~\eqref{eq:SDE}, $P_tf(x)=\EE[f(X_t^x)]$, and $P_t^*\delta_x$ for the law of $X_t^x$, that is, $P_t^*\delta_x(A)=\P(X_t^x\in A)$. The general notation is introduced in Subsection~\ref{subsec:notation}. It is well known that the generator of $(P_t)_{t\ge0}$,
\[(\mathcal{L}_p f)(z)=\frac{1}{2}f''(z)-f'(z)z^{[p-1]},\quad\text{for sufficiently smooth $f$},\]
possesses a spectral gap, see~\cite[Chapter~4]{BGL:14}.
Hence, the process $\{X_t^x\}_{t\geq 0}$ is uniquely ergodic with an invariant probability measure $\mu_p$ given by
\[
\mu_p(\ud z)=\frac{1}{\mathcal{Z}_p}\exp(-2\Psi_p(z))\ud z, \quad \textrm{ where }\quad \Psi_p(z)=\frac{1}{p}|z|^{p},\,\, z\in \mathbb{R},
\]
and $\mathcal{Z}_p$ is a normalizing constant.
The spectral gap follows from the symmetry of $\mathcal L_p$ in $L^2(\mu_p)$ and the Poincar\'e inequality for $\mu_p$, which holds since $\mu_p$ is log-concave~\cite{Bob:99}.
Moreover, the convergence also holds in the Wasserstein distance of order $r=2$, that is, for every initial datum $x\in\mathbb{R}$,
\begin{equation}
\lim\limits_{t\to \infty}\mathcal{W}_2(P_t^*\delta_x,\mu_p)=0,
\end{equation}
see Subsection~\ref{subsec:notation} for the definition of $\mathcal{W}_r$, $r\geq 1$.
The drift $z\mapsto-z^{[p-1]}$ is non-increasing. Hence two solutions driven by the same Brownian motion satisfy
\[
\frac{\ud}{\ud t}|X_t^x-X_t^y|^2=-2\left((X_t^x)^{[p-1]}-(X_t^y)^{[p-1]}\right)(X_t^x-X_t^y)\le0,
\]
that is, $|X_t^x-X_t^y|\le|x-y|$ for all $t\ge0$, $\P$-a.s. As in the proof of Lemma~\ref{lem:monotone} below, this synchronous coupling together with the invariance of $\mu_p$ implies that $t\mapsto\mathcal W_2(P_t^*\delta_x,\mu_p)$ is non-increasing. The Markov property alone would not suffice for this conclusion, see Remark~\ref{rem:markov}.

Given a prescribed error $\e>0$ we define the $\e$-mixing time with respect to $\mathcal{W}_2$ as
\begin{equation}\label{eq:mixing}
\tau_{\mathsf{mix}}(\e;x):=\inf\{t\geq 0: \mathcal{W}_2(P_t^*\delta_x,\mu_p)\leq \e\}.
\end{equation}
For $p\ge2$,\footnote{For $1<p<2$ the drift is not strongly dissipative at infinity, and a bound of the form~\eqref{eq:coupling} with a polynomial factor does not hold, see Appendix~B.} using~\cite{Eberle,K:18} one can show the existence of constants $C_p>0$ and $\lambda_p>0$ such that
\begin{equation}\label{eq:coupling}
\mathcal{W}_2(P_t^*\delta_x,P_t^*\delta_y)\leq C_pe^{-\lambda_p t}K_p(x,y)\quad \textrm{ for all }\quad x,y\in \mathbb{R},
\end{equation}
where $K_p(x,y):=1+|x|^p+|y|^{p}$ for $x,y\in \mathbb{R}$.
Combining~\eqref{eq:coupling} with the disintegration inequality of Lemma~\ref{lem:disi} for $r=2$ and with the inversion of the rate in Lemma~\ref{lem:inversion}, one obtains, for all $x\in\mathbb R$ and $\e\in(0,C_p\bar K_p(x))$ with $\bar K_p(x):=(\int_{\mathbb R}K_p(x,y)^2\,\mu_p(\ud y))^{1/2}$,
\begin{align*}
\tau_{\mathsf{mix}}(\e;x)\leq
\frac{1}{\lambda_p}\log\left(\frac{1}{\e}\right)+
\frac{1}{\lambda_p}\log\left(C_p \bar K_p(x)\right).
\end{align*}
The details are given in Appendix~B.
On the other hand, the deterministic equation
\[\frac{\ud}{\ud t} u_t^x=-(u_t^x)^{[p-1]},\quad u_0^x=x\in\R,\]
has, for $x\ne0$, the unique solution
\[
u_t^x=\begin{cases}
x(1+(p-2)|x|^{p-2}t)^{-1/(p-2)},&p>2,\\
xe^{-t},&p=2,\\
x\bigl(1-(2-p)|x|^{p-2}t\bigr)_+^{1/(2-p)},&1<p<2,
\end{cases}
\]
where $a_+:=\max\{a,0\}$. For $x=0$ the solution is identically zero. The third formula, valid only for $1<p<2$, includes the continuation by zero after the solution reaches zero.
The solution has the polynomial decay $u_t^x\to 0$ as $t\to\infty$ with optimal rate $t^{-\frac{1}{p-2}}$ for $p>2$. If
$1<p<2$, the solution reaches zero at the finite time $t_\ast=\frac{|x|^{2-p}}{2-p}$ and stays zero afterwards.
For $p=2$, the decay is exponential.
Thus, for $p>2$, the deterministic solution decays only polynomially, whereas the law of the stochastic solution converges exponentially fast to $\mu_p$. This is an example of \emph{stabilization by noise}, see~\cite{ACW:83} for linear systems.

\subsection{The stochastic $p$-Laplace equation}\label{subsec:pLap}

Let $U\subset\R^d$ be a possibly unbounded Poincar\'e domain\footnote{Examples of Poincar\'e domains are domains with finite width, that is, domains that fit between two parallel $(d-1)$-dimensional hyperplanes.}, that is, a non-empty domain on which Poincar\'e's inequality holds for the elements of $W^{1,p}_0(U)$. Denote $H:=L^2(U)$. Let $\{W_t\}_{t\ge 0}$ be a cylindrical Wiener process on a separable Hilbert space $H_0$, modeled on a filtered probability space $(\Omega,\Fcal,\{\Fcal_t\}_{t\ge 0},\P)$ satisfying the usual conditions. Let $B\in L_2(H_0,H)$, where $L_2(H_0,H)$ denotes the space of Hilbert--Schmidt operators from $H_0$ to $H$. Let $T>0$. The parabolic stochastic $p$-Laplace equation, $p>1$, $t\in [0,T]$, with zero Dirichlet boundary conditions and additive noise is given by
\begin{equation}\label{eq:pLaplace}
\ud X_t=\nabla\cdot[|\nabla X_t|^{p-2}\nabla X_t]\,\ud t+B\,\ud W_t,\quad X_0=x_0\in L^2(U),
\end{equation}
with $X_t\vert_{\partial U}=0$ for $t>0$.
The initial datum could be chosen in $L^1(U)$, see~\cite{SZ:21}, but we will only work with $L^2$-initial data for convenience.

\smallskip\noindent\textit{Assumptions.} We assume our domains to be Poincar\'e domains only for well-posedness. All quantitative statements in Sections~\ref{sec:overview} and~\ref{sec:mixing} are made in the following fixed setting. Where a cited result needs more or is formulated in a different setting, this is stated at the corresponding place.
\begin{enumerate}[(S1)]
\item $U\subset\R^d$ is a bounded domain with sufficiently smooth boundary, and the boundary condition is homogeneous Dirichlet. Convexity of $U$ is additionally assumed where we use~\cite{GT:16-1,GT:16,SST:23}. The results of~\cite{Liu:09-1,W:15-1,W:15} on non-degenerate noise are restricted to $d=1$.
\item $H=L^2(U)$ with norm $\|\cdot\|_H$, and $V_p:=W^{1,p}_0(U)\cap H$, which we use only for $p\ge\frac{2d}{d+2}$, with norm $v\mapsto\|\nabla v\|_{L^p(U)}$. All distances between probability measures are taken with respect to $\|\cdot\|_H$.
\item The noise is additive and the operator $B$ is fixed. The only exceptions are the remarks on~\cite{GT:16-1} and~\cite{W:15-1} in Section~\ref{sec:overview}, which also allow multiplicative noise.
\end{enumerate}
We use the following assumptions on $B$. Each of (N1)--(N3) implies (N0).
\begin{enumerate}[(N1)]
\setcounter{enumi}{-1}
\item $B\in L_2(H_0,H)$. This allows degenerate noise. For instance, $B$ may have a non-trivial kernel or finite rank, or vanish.
\item $B\in L_2(H_0,W^{1,2}_0(U))$, that is, $B$ maps $H_0$ into $W^{1,2}_0(U)$ and is Hilbert--Schmidt with respect to the norm of $W^{1,2}_0(U)$.
\item $\sum_{k}\|Be_k\|^2_{L^\infty(U)}<\infty$ for some orthonormal basis $\{e_k\}$ of $H_0$.
\item $d=1$, $U=(0,l)$, $H_0=H$, and $Be_i=q_ie_i$ for the eigenbasis $e_i=\sqrt{2/l}\,\sin(i\pi\cdot/l)$ of the Dirichlet Laplacian, where $q_i>0$, $\sum_iq_i^2<\infty$ and $q_i\ge\sigma i^{-1}$ for some $\sigma>0$ and all $i$. This is the non-degeneracy assumption of~\cite{W:15}, where, owing to a different normalization of the noise, $\sigma$ is replaced by $\sigma/\sqrt l$. For $l=1$ it is contained in the assumption of~\cite[Proposition 6.2]{W:15-1}, and it covers $B=(-\Delta)^{-\theta}$, $\theta\in(\frac14,\frac12]$, as in~\cite[Example 3.3]{Liu:09-1}.
\end{enumerate}
Several arguments, such as the pathwise contraction in Lemma~\ref{lem:monotone}, use the additive noise assumption in an essential way.

\subsection{Known results and open problems}

It has been proved by Gess and the second author~\cite{GT:14,GT:16} and by Liu~\cite{Liu:09-1} that the stochastic $p$-Laplace equation with degenerate noise is ergodic in all space dimensions and for every $p\in [1,\infty)$, including the borderline case of the multi-valued total variation flow $p=1$. For $p<2$ this requires additional assumptions on $B$ and, in~\cite{GT:16}, Neumann boundary conditions, see Subsection~\ref{subsec:degnoise}.

Known quantitative results are discussed in Section~\ref{sec:overview}, where the quantitative rates for the singular stochastic $p$-Laplace equation, that is, $p\in (1,2)$, depend on the spatial dimension. It is open how these dimensional restrictions can be relaxed. A striking result is stabilization by noise for $p>2$. In one spatial dimension and for non-degenerate noise, the stochastic $p$-Laplace equation is exponentially ergodic in total variation, uniformly in the initial datum. This was proved by Liu~\cite[Theorem 1.5(ii) and Example 3.3]{Liu:09-1} by means of dimension-free Harnack inequalities, and extended by Wang, in~\cite{W:15-1} by asymptotic couplings by reflection and in~\cite{W:15} with explicit rates. In contrast, the deterministic $p$-Laplace equation with $p>2$ converges at a polynomial rate.

The main open question concerns the singular case $1<p<2$ with degenerate additive noise. With homogeneous Dirichlet data, the solutions of the deterministic singular equation with $1<p<2$ and $p\ge\frac{2d}{d+2}$ vanish after a finite time, so that its mixing times stay bounded as $\eps\downarrow0$, see~\cite[Chapter~VII, Section~2]{DB:93} and Proposition~\ref{prop:extinction}. For $1<p<2$, $p>\frac{2d}{d+2}$, and degenerate noise, we show in Proposition~\ref{prop:singular} that the stochastic mixing times in every $\mathcal W_r$ are of order at most $\eps^{-\delta}$ for every $\delta>0$. It is open whether they are of logarithmic order, as for $p=2$. A list of open problems is given in Section~\ref{sec:conclusion}.

\subsection{Contributions and organization of the paper}\label{subsec:contrib}

\smallskip\noindent\textit{Purpose and method.} The purpose of this paper is to provide a recipe for deriving bounds on mixing times from quantitative convergence estimates, and to apply it systematically to the $p$-Laplace equation. The recipe has four ingredients.
\begin{enumerate}[(i)]
\item A contraction estimate for two solutions driven by the same noise, that is, for the synchronous coupling, see Lemma~\ref{lem:monotone}.
\item Moment bounds for the invariant measure, see Lemma~\ref{lem:moments}.
\item The disintegration inequality of Lemma~\ref{lem:disi}. It turns estimates for two solutions with different initial data into estimates for the distance to equilibrium.
\item The inversion of the resulting rate into a bound on the mixing time, see Lemma~\ref{lem:inversion}.
\end{enumerate}

\smallskip\noindent\textit{Nature of the results.} Since a large part of this paper is an overview, we distinguish three kinds of statements.
\begin{itemize}
\item \emph{Known results.} Section~\ref{sec:overview}, including Table~\ref{tab:rates}, collects convergence estimates from the literature together with their hypotheses. Apart from Lemma~\ref{lem:moments} and the short derivations of~\eqref{eq:SSTmu} and~\eqref{eq:expW1} from~\cite{SST:23} and~\cite{W:15}, no new estimates are proved there. Lemma~\ref{lem:moments} extends the moment bound of~\cite[Proposition~2.1]{BDP:06} by the same argument.
\item \emph{Direct consequences.} The mixing-time bounds in Table~\ref{tab:mix} whose source is an equation number follow from known estimates by Lemma~\ref{lem:inversion}. They require no new estimate. The same holds for the deterministic singular case in Proposition~\ref{prop:extinction}, which follows from~\cite{DB:93,BB:20}. The point is to identify the metric, the class of initial data and the moment assumptions under which the inversion is valid.
\item \emph{New statements.} The following results are proved here.
\begin{enumerate}[(1)]
\item The monotonicity of $t\mapsto\mathcal W_r(P_t^*\delta_x,\mu_*)$ for every $r\ge1$, which Lemma~\ref{lem:monotone} derives from the pathwise contraction.
\item For $p=2$, Proposition~\ref{prop:heat} gives bounds without smallness or spatial regularity assumptions on $B$, and Proposition~\ref{prop:heatlower} gives a lower bound of the same order $\lambda_1^{-1}\log(1/\eps)$ for every $x$ with $\langle x,e_1\rangle_H\neq0$.
\item For $p>2$, Proposition~\ref{prop:degenerate} gives the bound~\eqref{eq:Liu0} with an explicit constant, and Remark~\ref{rem:Wr} gives exponential mixing in every $\mathcal W_r$ for non-degenerate noise in one dimension.
\item For $1<p<2$, $p>\frac{2d}{d+2}$, and degenerate noise, Proposition~\ref{prop:singular} gives convergence in every $\mathcal W_r$, $r\ge1$, faster than any power of $t$, and hence $\tau_{\mathsf{mix},\mathcal W_r}(\eps;x)=O(\eps^{-\delta})$ for every $\delta>0$. This removes the restrictions $p\ge\sqrt2$ and $r=1$ of the bound that follows from~\cite{LT:11}, and improves its rate $t^{-1/2}$.
\end{enumerate}
In Table~\ref{tab:mix}, the mixing-time bounds among these results are referenced by the corresponding proposition or remark.
\end{itemize}

\smallskip\noindent\textit{Relation to~\cite{BT:24}.} In~\cite{BT:24} we studied ergodicity and mixing in $\mathcal W_2$ for locally monotone stochastic evolution equations with additive noise, under a dissipativity condition combined with a smallness condition on the noise. The resulting mixing times are logarithmic in $\eps$. The $p$-Laplace equation with $p\neq2$ is outside the scope of~\cite{BT:24}, and only the linear case $p=2$ is covered there. Compared with~\cite{BT:24}, the present paper differs in the following respects.
\begin{enumerate}[(a)]
\item It treats $p>1$ within the dimensional ranges of Table~\ref{tab:mix}, including regimes with polynomial rates or with solutions that vanish after a finite time.
\item It works with several distances, namely $\mathcal W_r$ for every $r\ge1$, a H\"older-type transport distance $D_\gamma$, and total variation, instead of $\mathcal W_2$ only.
\item For $p=2$ on bounded domains, it removes the smallness and spatial regularity assumptions on $B$ from the bound of~\cite[Theorem 2.3]{BT:24}.
\item In the singular regime, it compares the bounded mixing times of the deterministic equation, whose solutions vanish after a finite time, with stochastic upper bounds of order $\eps^{-\delta}$ for every $\delta>0$.
\end{enumerate}
We do not claim sharper constants than in the cited sources for the known results, see Remark~\ref{rem:constants}.

\smallskip\noindent\textit{Organization.} Section~\ref{sec:overview} recalls the solution concepts, fixes the notation for semigroups, distances and mixing times, proves moment bounds of all orders for the invariant measure in Lemma~\ref{lem:moments}, and collects the known convergence rates in Table~\ref{tab:rates}. Section~\ref{sec:mixing} proves monotonicity, the inversion lemma and the mixing-time bounds, which are summarized in Table~\ref{tab:mix}. Section~\ref{sec:conclusion} lists open problems. Appendix~A states the disintegration inequality, and Appendix~B gives the details of the example in Subsection~\ref{subsec:motivation}.

\section{Overview of known results}\label{sec:overview}

Let us collect known results on the stochastic and deterministic parabolic $p$-Laplace equation, see also~\cite{CGT:25}.
Throughout, the assumptions (S1)--(S3) are in force.

Let us define what we mean by a solution to~\eqref{eq:pLaplace}. For $p\ge 2$, that is, the degenerate $p$-Laplace case, we can
work with the following pathwise notion of a solution.
\begin{definition}\label{def:var}
Let $p\ge 2$. We say that $\{X_t\}_{t\in [0,T]}$ is a \emph{variational solution} to~\eqref{eq:pLaplace} for initial data $x_0\in H$ if
\begin{enumerate}[(i)]
\item $X\in L^2(\Omega;C([0,T];L^2(U)))\cap L^p([0,T]\times\Omega;W^{1,p}_0(U))$,
\item for any $v\in W^{1,p}_0(U)$, $t\in [0,T]$, $\P$-a.s.
\[\int_U X_t v\,\ud z=\int_U x_0 v\,\ud z-\int_0^t\int_U |\nabla X_s|^{p-2}(\nabla X_s \cdot \nabla v)\,\ud z\, \ud s +\int_U v B \, W_t \,\ud z.\]
\end{enumerate}
\end{definition}
Under (N0), the existence and uniqueness of variational solutions in all space dimensions $d\ge 1$ has been proved in~\cite{Liu:09,PR:07}. Definition~\ref{def:var} extends without change to $1<p<2$ with $p\ge\frac{2d}{d+2}$, where $V_p\subset H\subset V_p^*$ is a Gelfand triple, and for these $p$ variational solutions exist and are unique by~\cite[Section~1]{Liu:09}. For $p\ge\frac{2d}{d+2}$ we always use variational solutions. The SVI solutions of Definition~\ref{def:svi} below are needed only for $1<p<\frac{2d}{d+2}$.

For $p\in (1,2)$, that is, the singular $p$-Laplace case, we can
work with the following notion of a weak variational solution.
\begin{definition}\label{def:svi}
Let $p\in (1,2)$. We say that $\{X_t\}_{t\in [0,T]}$ is a \emph{stochastic variational inequality (SVI) solution} to~\eqref{eq:pLaplace} for initial data $x_0\in L^2(U)$ if
\begin{enumerate}[(i)]
\item $X\in L^2([0,T]\times\Omega;L^2(U))$ such that
\[\operatorname{ess\,sup}\displaylimits_{t\in [0,T]}\EE\left[\|X_t\|^2_{L^2(U)}\right]+\EE\int_0^T\int_U|\nabla X_s|^p\,\ud z\,\ud s<\infty,\]
\item for any $z_0\in L^2(U)$, every $Z\in L^2([0,T]\times\Omega;W^{1,2}_0(U))$, and every progressively measurable $G\in L^2([0,T]\times\Omega;L^2(U))$ such that
\[Z_t=z_0+\int_0^t G_s\,\ud s+B\,W_t,\]
we have that for a.e. $t\in [0,T]$
\begin{align*}&\EE\left[\|X_t-Z_t\|^2_{L^2(U)}\right]+\frac{2}{p}\EE\int_0^t\int_U|\nabla X_s|^p\,\ud z\,\ud s\\
\le &\EE\left[\|x_0-z_0\|^2_{L^2(U)}\right]+\frac{2}{p}\EE\int_0^t\int_U|\nabla Z_s|^p\,\ud z\,\ud s\\
&-2\EE\int_0^t\int_U G_s(X_s-Z_s)\,\ud z\,\ud s.\end{align*}
\end{enumerate}
\end{definition}
The test processes $Z$ in Definition~\ref{def:svi} satisfy $Z_t=z_0+\int_0^tG_s\,\ud s+BW_t$ and are required to take values in $W^{1,2}_0(U)$. A rich class of such processes exists only if $BW_t\in W^{1,2}_0(U)$. Therefore this notion is used under (N1). Under the corresponding assumption $B\in L_2(H_0,W^{1,2}(U))$, existence and uniqueness of SVI solutions for $p\in(1,2)$ in all space dimensions $d\ge1$ was first proved in~\cite{GT:16-1} on bounded convex smooth domains $U$, there for homogeneous Neumann boundary conditions and possibly multiplicative noise. For homogeneous Dirichlet boundary conditions on bounded convex domains with piecewise $C^2$-boundary and under (N1), existence and uniqueness of generalized solutions follow from~\cite[Theorem 4.6]{GT:14}, see~\cite[Definition 2.4, Theorem 2.5 and Corollary 2.7]{SST:23}. In this case, by the solution of~\eqref{eq:pLaplace} we always mean this generalized solution. See~\cite{GT:14} for the notion of \emph{limit solutions}, which are also known as \emph{generalized solutions}, or as \emph{SOLA}\footnote{\emph{SOLA} is an abbreviation for \emph{solutions obtained as limits of approximations.}}. Limit solutions are always SVI solutions, see~\cite[Appendix C]{GT:14}.

\subsection{Semigroups, distances and mixing times}\label{subsec:notation}

The transition semigroup $(P_t)_{t\ge 0}$ associated with~\eqref{eq:pLaplace} is Markovian and Feller. It acts as follows
\[
P_t F(x):=\EE[F(X_t^x)]\quad
\textrm{ for any }\quad F\in B_b(H)\quad \textrm{ and }\quad x\in H,
\]
where $B_b(H):=\{F:H\to \mathbb{R}:\,F\textrm{ is bounded and Borel measurable}\}$. See~\cite[Proposition~4.3.5]{LR:15} for a proof of the Markov property.
 For a semigroup $(P_t)_{t\ge 0}$, we define the dual semigroup $(P_t^\ast)_{t\ge 0}$ acting on $\Mcal_1(H,\Bcal(H)):=\{\mu:\Bcal(H)\to [0,1]:\,\mu \textrm{ is a probability measure}\}$ by
\[P_t^\ast \mu(A):=\int_H P_t \mathbbm{1}_A(x)\,\mu(\ud x)\quad \textrm{ for any }\quad A\in\Bcal(H),
\]
where $\Bcal(H)$ denotes the Borel sets of $H$ and $\mathbbm{1}_A$ denotes the indicator function of the set $A$, see~\cite{DPZ:96} for details.

A measure $\mu\in\Mcal_1(H,\Bcal(H))$ is said to be \emph{invariant} for the semigroup $(P_t)_{t\ge 0}$ if $P^\ast_t\mu=\mu$ for all $t\ge 0$.

We distinguish four objects. The solution $X_t^x$ is an $H$-valued random variable. Its law $P_t^*\delta_x$ is a probability measure on $H$ with $P_t^*\delta_x(A)=\P(X_t^x\in A)$. The semigroup $(P_t)_{t\ge0}$ acts on functions, and the dual semigroup $(P_t^*)_{t\ge0}$ acts on probability measures. Thus $P_tF(x)=\EE[F(X_t^x)]=\langle F,P_t^*\delta_x\rangle$, where $\langle F,\mu\rangle:=\int_HF\,\ud\mu$. Whenever existence and uniqueness of an invariant probability measure hold under the hypotheses in use, we denote it by $\mu_*$. When the dependence on the data matters, we write $P_t=P_t^{p,B,U}$ and $\mu_*=\mu_*^{p,B,U}$. In the deterministic case $B=0$ we have $P_t^*\delta_x=\delta_{u_t^x}$, where $u_t^x$ is the solution of~\eqref{eq:detpLaplace} below, and $\mu_*=\delta_0$.

\smallskip\noindent\textit{Distances.} For $r\ge1$ let $\mathcal P_r(H)$ denote the set of probability measures on $(H,\mathcal{B}(H))$ with finite $r$th absolute moment.
The Wasserstein distance of order $r$ between $\mu,\nu\in\mathcal P_r(H)$ is defined by
\begin{equation}\label{eq:defWr}
\mathcal{W}_r(\mu,\nu):=\left(\inf_{\Gamma\in \mathcal{C}(\mu,\nu)}
\int_{H\times H} \|x-y\|^r_H\,\Gamma(\ud x,\ud y)\right)^{1/r},
\end{equation}
where the infimum is taken over the set $\mathcal C(\mu,\nu)$ of all couplings $\Gamma$ between $\mu$ and $\nu$, that is,
for any $A\in \mathcal{B}(H)$
\begin{equation}\label{eq:defc}
\Gamma(A\times H)=\mu(A)\quad \textrm{ and }\quad \Gamma(H\times A)=\nu(A).
\end{equation}
In particular, $\mathcal W_r(\mu,\nu)\le(\int_{H\times H}\|x-y\|_H^r\,\Gamma(\ud x,\ud y))^{1/r}$ for any coupling $\Gamma$.
By~\cite[Theorem 6.9]{villani}, convergence in $\mathcal W_r$ on $\mathcal P_r(H)$ is equivalent to weak convergence together with convergence of the $r$th absolute moments. Weak convergence alone does not suffice. By Jensen's inequality, $\mathcal W_r(\cdot,\cdot)\le\mathcal W_{r'}(\cdot,\cdot)$ for $1\le r\le r'$.
For $0<\gamma\le1$ and probability measures $\mu,\nu$ with finite $\gamma$th moments, let
\[
D_\gamma(\mu,\nu):=\sup_{\|F\|_\gamma\le1}|\langle F,\mu-\nu\rangle|,\qquad
\|F\|_\gamma:=\sup_{x\not=y}\frac{|F(x)-F(y)|}{\|x-y\|_{H}^\gamma}.
\]
By Kantorovich--Rubinstein duality~\cite[Chapter 5]{villani}, $D_\gamma$ is the transport distance with cost $\|x-y\|_H^\gamma$. In particular, $D_1(\cdot,\cdot)=\mathcal W_1(\cdot,\cdot)$, and $D_\gamma(\cdot,\cdot)\le\mathcal W_1(\cdot,\cdot)^\gamma$ by Jensen's inequality. We write $\Lip(F):=\|F\|_1$ for the Lipschitz constant of $F$. For $\gamma<1$, $D_\gamma$ is not a Wasserstein distance of order $r\ge1$, and a rate in $D_\gamma$ does not yield the same rate in $\mathcal W_r$. Finally, $D_{\rm TV}(\mu,\nu):=\sup_{A\in\mathcal B(H)}|\mu(A)-\nu(A)|$ is the total variation distance. It satisfies $2D_{\rm TV}(\mu,\nu)=\sup_{\|F\|_\infty\le1}|\langle F,\mu-\nu\rangle|=|\mu-\nu|(H)$, where $|\mu-\nu|$ denotes the total variation measure of $\mu-\nu$. It requires no moments and is not comparable to $\mathcal W_r$ on $H$ in general.

\smallskip\noindent\textit{Mixing times.} Let $D\in\{\mathcal W_r,D_\gamma,D_{\rm TV}\}$, $\eps>0$ and $x\in H$. The $\eps$-mixing time of~\eqref{eq:pLaplace} with respect to $D$ is
\begin{equation}\label{eq:mixdef}
\tau_{\mathsf{mix},D}^{p,B,U}(\eps;x):=\inf\{t\ge0:D(P_t^*\delta_x,\mu_*^{p,B,U})\le\eps\},\qquad \inf\varnothing:=+\infty.
\end{equation}
We suppress $p,B,U$ when they are fixed. In the settings of Lemma~\ref{lem:monotone} and Remark~\ref{rem:markov}, the distance $D(P_t^*\delta_x,\mu_*)$ is non-increasing in $t$, so it does not exceed $\eps$ after the mixing time. Remark~\ref{rem:envelope} discusses a variant of~\eqref{eq:mixdef} for distances that are not known to be monotone.

\subsection{The deterministic case}
\label{sec:dcplaplace}

In the deterministic case, $P_t F(x)=F(u^x_t)$, $t\ge 0$, where $F\in C_b(H)$ is an observable.
Here $\{u_t^x\}$ is the analytically weak solution to
\begin{equation}\label{eq:detpLaplace}
\partial_t u_t^x=\nabla\cdot[|\nabla u_t^x|^{p-2}\nabla u_t^x],\quad u_0^x=x\in L^2(U),\quad u_t^x\vert_{\partial U}=0\quad\text{for}\quad t>0.
\end{equation}

Juutinen and Lindqvist~\cite{JL:09} study the decay of solutions with time-independent lateral boundary values towards the stationary solution. They obtain polynomial rates for $p>2$ and exponential rates in the singular case. With homogeneous Dirichlet data, $1<p<2$ and $p\ge\frac{2d}{d+2}$, the solutions vanish after a finite time, see~\cite[Chapter~VII, Section~2]{DB:93} and Proposition~\ref{prop:extinction}.

For $p>2$, they prove in~\cite[Proposition 3.1]{JL:09} that there exists a constant $C=C(p,d,U)>0$ such that for every $m\in (1,\infty)$
\begin{equation}\label{eq:detpol}\|u_t^x-u_t^y\|_{L^m(U)}\le\|x-y\|_{L^m(U)}\wedge C\left(\frac{m-1}{(m+p-2)^p}\right)^{\frac{1}{2-p}} t^{-\frac{1}{p-2}},\quad\text{for all}\quad t>0.\end{equation}

For the singular range $1<p<2$ with $p\ge\frac{2d}{d+2}$, they prove in~\cite[Proposition 3.4]{JL:09} that
there exists a constant $C=C(p,d,U)>0$ such that for every $m\in (1,\infty)$
\begin{equation}\label{eq:detexp}\|u_t^x-u_t^y\|_{L^m(U)}\le\|x-y\|_{L^m(U)}e^{-\frac{CL(t)}{m}t},\quad\text{for all}\quad t>0,\end{equation}
where
\[L(t)=\left(\operatorname{ess\,sup}\displaylimits_{0\le s\le t}\left(\|u_s^x\|_{W^{1,p}_0}+\|u_s^y\|_{W^{1,p}_0}\right)\right)^{p-2}.\]
The estimate requires $x,y\in V_p\cap L^m(U)$. Since~\eqref{eq:detpLaplace} is the gradient flow in $H$ of the convex functional $v\mapsto\frac1p\int_U|\nabla v|^p\,\ud z$, the energy is non-increasing along solutions~\cite{Br:73}. Hence $L(t)\ge(\|\nabla x\|_{L^p}+\|\nabla y\|_{L^p})^{p-2}>0$ for all $t$, which gives an exponential rate depending on the energies of $x$ and $y$. In the present setting, the solutions even vanish after a finite time, see Proposition~\ref{prop:extinction}.

Furthermore, DiBenedetto~\cite{DB:93}, DiBenedetto and Herrero~\cite{DBH:89,DBH:90}, and Porzio~\cite{P:09,P:11} have obtained the following decay estimates, which are also called ultracontractive estimates, for the deterministic parabolic $p$-Laplace equation with zero Dirichlet boundary conditions, for $\frac{2d}{d+m}<p<d$, $m\ge 1$, and for some $C=C(d,m,p,U)>0$,
\[\|u_t^x\|_{L^\infty(U)}\le C t^{-\frac{d}{d(p-2)+pm}}\|x\|_{L^m(U)}^{\frac{pm}{d(p-2)+pm}},\quad t>0,\]%
and in the degenerate case $p>2$, the universal bound for some $C>0$
\[\|u_t^x\|_{L^\infty(U)}\le C t^{-\frac{1}{p-2}},\quad t>0,\]
holds uniformly for all $x$.
Moreover, for the singular case, namely if either $m>1$ and $\frac{2d}{d+m}\le p<2$, or $m=1$ and $\frac{2d}{d+1}<p<2$, the solution vanishes after a finite time, with upper bounds for this time given in~\cite[Theorems 1.5 and 1.6]{P:11}.

Finally, in the case of the heat equation, they obtain for
$p=2$, and for some $c_i>0$, $i=1,2$, $C>0$,
\[\|u_t^x\|_{L^\infty(U)}\le Ct^{-c_1}e^{-c_2 t}\|x\|_{L^{m}(U)},\quad t>0.\]
Here the factor $C\|x\|_{L^m(U)}$ necessarily depends on $x$. Indeed, by linearity, no bound $\|u_t^x\|_{L^\infty(U)}\le\rho(t)$ with $\rho(t)\to0$ as $t\to\infty$ can hold uniformly in $x\in H$.

\subsection{The case of degenerate noise}\label{subsec:degnoise}

Let $B$ satisfy (N0). In particular, the noise may be spatially degenerate. Let $\{X_t^x\}_{t\ge 0}$ be the solution to~\eqref{eq:pLaplace} with $X_0=x$.
The existence and uniqueness of invariant measures for the stochastic $p$-Laplace equation have been discussed in~\cite{Liu:09,LT:11,GT:14,GT:16,Liu:09-1,W:15-1,W:15,BDP:06}, among others.
The estimate~\eqref{eq:Liu} of Liu~\cite{Liu:09-1} below establishes the uniqueness of an invariant measure for all spatial dimensions for $p>2$. The existence of an invariant measure follows from coercivity and the Krylov--Bogoliubov method, cf.~\cite{DPZ:96}.
In~\cite{GT:16}, under the additional regularity assumptions (N2), $H_0=H$ and $B\in L_2(H,H^3(U))$ symmetric, and for homogeneous Neumann boundary conditions, Gess and the second author have proved the existence of a unique invariant measure in all spatial dimensions for $p\in [1,2)$, including the singular borderline case of the stochastic total variation flow $p=1$. Partial results for the singular case in low dimensions have been proved before in~\cite{GT:14,Liu:09,LT:11}.

For $p>2$, Liu showed the following bounds in~\cite[Theorem 1.3]{Liu:09-1}, which remain true for degenerate noise, for some $C>0$,
\begin{equation}\label{eq:Liu0}\|X_t^x-X_t^y\|_{L^2(U)}^2\le\|x-y\|^2_{L^2(U)}\wedge \left(Ct^{-\frac{2}{p-2}}\right),\quad t>0,\quad\P\text{-a.s.},\end{equation}
and
\begin{equation}\label{eq:Liu}|P_t F(x)-\langle F,\mu_*\rangle|\le C\Lip(F)t^{-\frac{1}{p-2}},\quad t>0,\end{equation}
which recovers the rate of the deterministic case. A pathwise bound of this form for drifts $A$ with $\langle A(u)-A(v),u-v\rangle\le-c\|\nabla(u-v)\|_{L^p(U)}^p$, which include the $p$-Laplacian with $p\ge2$ by~\cite[Section~12, (I)]{Lin:17}, is already obtained in the proof of~\cite[Proposition~2.2]{BDP:06}. The bound~\eqref{eq:Liu0} is pathwise, uniform in $x,y\in H$ and holds for every $t>0$. A self-contained proof with an explicit constant is given in Proposition~\ref{prop:degenerate}.

For $1<p<2$ with $p>\frac{2d}{d+2}$, in~\cite[Theorem 1.2 and Proposition 3.2]{LT:11}, Liu and the second author proved the following rate for some $C>0$, under (N0) and for all $x,y\in H$,
\begin{equation}\label{eq:LT}\begin{split}&\EE\left[\|X_t^x-X_t^y\|_{L^2(U)}^{\frac{2p}{2-p}}\right]\\
\le& C\|x-y\|^{\frac{2p}{2-p}}_{L^2(U)}t^{-\frac{p}{2-p}}\left(1+t^{-1}\|x\|_{L^2(U)}^2+t^{-1}\|y\|_{L^2(U)}^2\right),\quad t>0,\end{split}\end{equation}
and, by~\cite[Theorem 1.4]{LT:11}, if $p\ge\sqrt{2}$,
\begin{align*}&|P_t F(x)-\langle F,\mu_*\rangle|\\
\le& C\Lip(F)(1+\|x\|_{L^2(U)})t^{-\frac{1}{2}}\left(1+t^{-\frac{2-p}{2p}}(1+\|x\|_{L^2(U)})^{\frac{2-p}{p}}\right),\quad t>0,\end{align*}
and if $p\le \sqrt{2}$, for any $\gamma\in \left(0,\frac{p^2}{2}\right]$,
\begin{align*}&|P_t F(x)-\langle F,\mu_*\rangle|\\
\le& C\|F\|_\gamma (1+\|x\|_{L^2(U)}^\gamma)t^{-\frac{\gamma}{2}}\left(1+t^{-\frac{(2-p)\gamma}{2p}}(1+\|x\|_{L^2(U)})^{\frac{(2-p)\gamma}{p}}\right),\quad t>0,\end{align*}
where $\|F\|_\gamma$ is the $\gamma$-H\"older seminorm defined in Subsection~\ref{subsec:notation}.
The restrictions $p\ge\sqrt2$ and $\gamma\le p^2/2$ arise because only the moment $\int_H\|\nabla y\|_{L^p}^p\,\mu_*(\ud y)<\infty$ is used in~\cite{LT:11} when integrating~\eqref{eq:LT} against $\mu_*$. With the moments of all orders of Lemma~\ref{lem:moments} below, both restrictions can be removed, see Proposition~\ref{prop:singular}.

Using second-order estimates, Seib, Stannat and the second author refined this result in~\cite[Proposition~3.3, Theorem~5.1 and Corollary~5.3]{SST:23}. They assume (N1), a convex bounded domain $U$ with piecewise $C^2$-boundary, $1<p<2$, $p\ge2-\frac4d$, and initial data $x,y\in W^{1,2}_0(U)$, and obtain, for some $C_i>0$, $i=1,2,3$, and every $0<\gamma\le\gamma_*:=1\wedge\frac{8-2(2-p)}{(2-p)(4-p)(1\vee\frac{d}{2})}$, the estimate
\begin{equation}\label{eq:SST}\begin{split}&\EE\left[\|X_t^x-X_t^y\|^\gamma_{L^2(U)}\right]\\
\le& C_1 t^{-\frac{p}{2-p}}\|x-y\|_{L^2(U)}^\gamma\\
&\times\EE\left[\frac{1}{t}\int_0^t\left(C_2+C_3\left(\|X_s^x\|_{W^{1,2}_0(U)}^p+\|X_s^y\|_{W^{1,2}_0(U)}^p\right)\right)\,\ud s\right],\quad t>0.\end{split}\end{equation}
The estimate~\eqref{eq:SST} is contained in the proof of~\cite[Theorem 5.1]{SST:23}. Its rate does not depend on $\gamma$. By~\cite[Theorem 5.1 and Corollary 5.3]{SST:23}, it implies
\begin{align*}
&\limsup_{t\to\infty}\left(t^{\frac{p}{2-p}}\frac{|P_t F(x)-P_t F(y)|}{\|x-y\|^\gamma_{L^2(U)}}\right)\\
&\le C_1\|F\|_\gamma\left(C_2+2C_3\int_{L^2(U)}\|z\|^p_{W^{1,2}_0(U)}\,\mu_*(\ud z)\right),
\end{align*}
for any $0<\gamma\le\gamma_*$, with the constants of~\eqref{eq:SST}.
Assuming higher regularity of the noise and the initial datum, the dimensional restriction and the rate are improved. This result also improves previous results in the deterministic case regarding the dimensional restriction. See~\cite{FPS:23,FPS:22} for regularity theory of the deterministic $p$-Laplace equation that employs second-order estimates.
The estimate~\eqref{eq:SST} bounds the $\gamma$th moment of $\|X_t^x-X_t^y\|_{L^2(U)}$. It therefore concerns the distance $D_\gamma$, $0<\gamma\le\gamma_*$, and not $\mathcal W_r$, $r\ge1$, except $D_1=\mathcal W_1$ when $\gamma=1$. We now derive a rate to equilibrium from it. For $x\in W^{1,2}_0(U)$ let
\begin{multline*}
R_\gamma(x):=\sup_{t\ge1}\int_H\|x-y\|_H^\gamma\\
\times\EE\left[\frac1t\int_0^t\left(C_2+C_3\left(\|X_s^x\|_{W^{1,2}_0(U)}^p+\|X_s^y\|_{W^{1,2}_0(U)}^p\right)\right)\ud s\right]\mu_*(\ud y),
\end{multline*}
and assume that $R_\gamma(x)<\infty$. By invariance, $\mu_*=\int_HP_t^*\delta_y\,\mu_*(\ud y)$. Since $D_\gamma$ is a supremum of linear functionals, it follows that
\[
D_\gamma(P_t^*\delta_x,\mu_*)\le\int_HD_\gamma(P_t^*\delta_x,P_t^*\delta_y)\,\mu_*(\ud y).
\]
For fixed $y$, the pair $(X_t^x,X_t^y)$ is a coupling of $P_t^*\delta_x$ and $P_t^*\delta_y$. As $D_\gamma$ is the transport distance with cost $\|a-b\|_H^\gamma$, we get $D_\gamma(P_t^*\delta_x,P_t^*\delta_y)\le\EE\|X_t^x-X_t^y\|_H^\gamma$. Under (N1), $\mu_*(W^{1,2}_0(U))=1$ by~\cite[Corollary 4.4]{SST:23}. Hence~\eqref{eq:SST} applies for $\mu_*$-a.e.\ $y$, and integrating it against $\mu_*(\ud y)$ gives
\begin{equation}\label{eq:SSTmu}
D_\gamma(P_t^*\delta_x,\mu_*)\le C_1R_\gamma(x)\,t^{-p/(2-p)},\qquad t\ge1.
\end{equation}
The condition $R_\gamma(x)<\infty$ is part of the corresponding table entries. The limit superior above holds for fixed $x$ and $y$, whereas integrating in $y$ requires a bound that is uniform in $t\ge1$, and $R_\gamma(x)<\infty$ provides it. For $B=0$, the same argument applies with $\mu_*=\delta_0$ and $x\in W^{1,2}_0(U)$ with $R_\gamma(x)<\infty$.

For $p>2$, and under the additional regularity assumption on $B$ discussed in~\cite[Remark 4.13]{G:12}, where, in particular, it is assumed that there exists a complete orthonormal system $\{e_k\}$ of $H_0$ such that
\[\sum_{k=1}^\infty \|Be_k\|^p_{W^{1,p}_0(U)}<\infty,\]
Gess~\cite{G:12} proved that
\[\sup_{s\in[0,T]}\EE\left[\|X_s^x\|^p_{W^{1,p}_0(U)}\right]+\EE\int_0^T\|\nabla\cdot(|\nabla X_s^x|^{p-2}\nabla X_s^x)\|^2_{L^2(U)}\,\ud s<\infty,\]
for any initial datum $x\in W^{1,p}_0(U)\cap L^2(U)$.

For $1<p<2$ with $p\ge2-\frac4d$ and assuming (N1), under convexity and boundedness of $U$, Seib, Stannat and the second author proved in~\cite[Corollary 4.4]{SST:23} that the unique invariant measure is concentrated on $W^{2,\frac{dp}{d+p-2}}_0(U)\cap W^{1,2}_0(U)$ and has $p$th moments in $W^{1,2}_0(U)$.

For $p\in[1,2)$ on bounded convex smooth domains $U$, assuming (N2), $H_0=H$ and that $B\in L_2(H,H^3(U))$ is symmetric,~\cite[Lemma 9]{GT:16} gives for any $d\ge 1$ and any $q\ge 2$ constants $C,c>0$ such that\footnote{The results in~\cite{GT:16} are formulated for homogeneous Neumann boundary conditions and data with zero mean. The authors of~\cite{GT:16} expect the Dirichlet case to be analogous.}
\[\EE\left[\|X_t^x\|^q_{L^q(U)}\right]+cq\EE\int_0^t\|X_s^x\|^{p+q-2}_{L^{p+q-2}(U)}\,\ud s\le \|x\|^q_{L^q(U)}+Cqt,\quad t>0.\]

\subsection{Moments of the invariant measure}\label{subsec:moments}

The following lemma provides moments of all orders for the invariant measure. It extends~\cite[Proposition~2.1]{BDP:06}, which treats $p\ge2$ and gives $\int_H(\|y\|_H^2+\|\nabla y\|_{L^p}^p)\,\mu_*(\ud y)<\infty$, to moments of all orders and to $\frac{2d}{d+2}<p<2$. The proof follows the same route.

\begin{lemma}[Moments of the invariant measure]\label{lem:moments}
Let $p>1$ with $p>\frac{2d}{d+2}$, assume (N0), and let $\mu_*$ be the unique invariant probability measure of~\eqref{eq:pLaplace}. Let $C_S>0$ be such that $\|v\|_H\le C_S\|\nabla v\|_{L^p(U)}$ for $v\in V_p$. Then
\begin{equation}\label{eq:moment}
\int_H\|\nabla y\|^p_{L^p(U)}\,\mu_*(\ud y)\le\frac12\|B\|^2_{L_2(H_0,H)},\qquad
\int_H\|y\|^p_H\,\mu_*(\ud y)\le\frac{C_S^p}{2}\|B\|^2_{L_2(H_0,H)},
\end{equation}
and for every $q\ge2$ there is a constant $K_q$, depending only on $q$, $p$, $C_S$ and $\|B\|_{L_2(H_0,H)}$, such that
\begin{equation}\label{eq:allmoments}
\int_H\|y\|_H^{q-2+p}\,\mu_*(\ud y)\le K_q.
\end{equation}
In particular, $\mu_*\in\mathcal P_r(H)$ for every $r\ge1$, and $P_t^*\delta_x\in\mathcal P_r(H)$ for every $x\in H$, $t\ge0$ and $r\ge1$.
\end{lemma}

The constant $C_S$ exists by the Sobolev embedding for $1<p<2$ with $p>\frac{2d}{d+2}$, and by Poincar\'e's and H\"older's inequalities on the bounded domain $U$ for $p\ge2$. Existence and uniqueness of $\mu_*$ in this range of $p$ follow from~\cite{Liu:09-1} for $p>2$, from~\cite{LT:11} for $p<2$, and from~\cite[Chapter 6]{DPZ:96} for $p=2$.

\begin{proof}
\textit{Step 1. Energy identity.} For these $p$, $V_p\subset H\subset V_p^*$ is a Gelfand triple, and the solution of~\eqref{eq:pLaplace} is a variational solution. It\^o's formula~\cite[Theorem 4.2.5]{LR:15} for $\|X_t^x\|_H^2$ gives
\[
\|X_t^x\|_H^2=\|x\|_H^2-2\int_0^t \|\nabla X_s^x\|_{L^p}^p\,\ud s+t\|B\|^2_{L_2(H_0,H)}+N_t,
\]
where $N_t:=2\int_0^t\langle X_s^x,B\,\ud W_s\rangle_H$ is a local martingale with quadratic variation $\langle N\rangle_t=4\int_0^t\|B^*X_s^x\|_{H_0}^2\,\ud s$.

\textit{Step 2. Higher powers.} Fix $q\ge2$ and $\delta>0$, and let $S_t:=\|X_t^x\|_H^2$. We apply It\^o's formula to $(\delta+S_t)^{q/2}$. By Step~1, the drift of $S$ is $-2\|\nabla X_t^x\|_{L^p}^p+\|B\|^2_{L_2(H_0,H)}$, and the density of its quadratic variation is $4\|B^*X_t^x\|_{H_0}^2\le4\|B\|^2_{L_2(H_0,H)}S_t$. Since $S_t\le\delta+S_t$, the second order term satisfies
\[
\frac12\cdot\frac q2\Big(\frac q2-1\Big)(\delta+S_t)^{\frac q2-2}\cdot4\|B^*X_t^x\|_{H_0}^2\le\frac q2(q-2)\|B\|^2_{L_2(H_0,H)}(\delta+S_t)^{\frac q2-1}.
\]
Hence, in the integrated sense and with a local martingale $N^{(q,\delta)}$,
\[
\ud(\delta+S_t)^{q/2}\le\frac q2(\delta+S_t)^{\frac q2-1}\Big(-2\|\nabla X_t^x\|_{L^p}^p+(q-1)\|B\|^2_{L_2(H_0,H)}\Big)\ud t+\ud N^{(q,\delta)}_t .
\]
Letting $\delta\downarrow0$, we obtain
\[
\ud\|X_t^x\|_H^q\le\left(-q\|X_t^x\|_H^{q-2}\|\nabla X_t^x\|_{L^p}^p+\frac{q(q-1)}{2}\|B\|^2_{L_2(H_0,H)}\|X_t^x\|_H^{q-2}\right)\ud t+\ud N^{(q)}_t
\]
with a local martingale $N^{(q)}$.

\textit{Step 3. Moment bound along the solution.} By the choice of $C_S$, $\|\nabla X\|_{L^p}^p\ge C_S^{-p}\|X\|_H^p$, so the first term of the drift is at most $-qC_S^{-p}\|X_t^x\|_H^{q-2+p}$. For the second term we use Young's inequality. It gives a constant $K'_q$, depending only on $q$, $p$, $C_S$ and $\|B\|_{L_2(H_0,H)}$, such that
\[
\frac{q(q-1)}{2}\|B\|^2_{L_2(H_0,H)}s^{q-2}\le\frac q2C_S^{-p}s^{q-2+p}+K'_q,\qquad s\ge0.
\]
Hence the drift is at most $-\frac q2C_S^{-p}\|X_t^x\|_H^{q-2+p}+K'_q$. We stop the process at $\sigma_n:=\inf\{t\ge0:\|X_t^x\|_H\ge n\}$ and take expectations, which removes the stochastic integral. Since $X^x$ has continuous paths in $H$, $\sigma_n\to\infty$ as $n\to\infty$. Fatou's lemma and monotone convergence give
\begin{equation}\label{eq:itomoment}
\EE\|X_t^x\|_H^q+\frac q2C_S^{-p}\,\EE\int_0^t\|X_s^x\|_H^{q-2+p}\,\ud s\le\|x\|_H^q+K'_qt,\qquad t\ge0.
\end{equation}
The same argument applied to the identity in Step~1 gives
\[
\EE\|X_t^x\|_H^2+2\EE\int_0^t\|\nabla X_s^x\|_{L^p}^p\,\ud s\le\|x\|_H^2+t\|B\|^2_{L_2(H_0,H)} .
\]

\textit{Step 4. Time averages.} For $T>0$ let $\nu_T:=T^{-1}\int_0^TP_s^*\delta_0\,\ud s$. We apply Step~3 with $x=0$, divide by $T$ and drop the first terms on the left-hand sides. This gives, for all $T>0$,
\[
\int_H\|\nabla y\|_{L^p}^p\,\nu_T(\ud y)\le\frac12\|B\|^2_{L_2(H_0,H)},\qquad\int_H\|y\|_H^{q-2+p}\,\nu_T(\ud y)\le\frac{2C_S^pK'_q}{q}.
\]
Since $p>\frac{2d}{d+2}$ and $U$ is bounded, $V_p$ embeds compactly into $H$. Hence the sets $\{y\in V_p:\|\nabla y\|_{L^p}\le R\}$ are compact in $H$. By Markov's inequality, they have $\nu_T$-measure at least $1-\|B\|^2_{L_2(H_0,H)}/(2R^p)$ for all $T>0$, so the family $(\nu_T)_{T>0}$ is tight. By the Feller property, every weak limit point of $(\nu_{T_n})$ with $T_n\to\infty$ is invariant~\cite[Theorem 3.1.1]{DPZ:96}, hence equal to $\mu_*$.

\textit{Step 5. Passage to the limit.} The functional $y\mapsto\|\nabla y\|_{L^p}^p$, set to $+\infty$ for $y\in H\setminus V_p$, is lower semicontinuous on $H$. Indeed, bounded sets of $V_p$ are weakly relatively compact, and the norm of $V_p$ is weakly lower semicontinuous. The functional $y\mapsto\|y\|_H^{q-2+p}$ is continuous on $H$. For a non-negative lower semicontinuous function $g$ on $H$, weak convergence $\nu_{T_n}\to\mu_*$ implies $\int_Hg\,\ud\mu_*\le\liminf_n\int_Hg\,\ud\nu_{T_n}$. Hence the bounds of Step~4 pass to $\mu_*$. This proves the first inequality in~\eqref{eq:moment}, and~\eqref{eq:allmoments} with $K_q:=2C_S^pK'_q/q$. The second inequality in~\eqref{eq:moment} follows from the first one and $\|y\|_H\le C_S\|\nabla y\|_{L^p}$.

\textit{Step 6. Conclusion.} Since $q\ge2$ is arbitrary, $\mu_*$ has finite moments of every order $r\ge p$, and hence of every order $r\ge1$. By~\eqref{eq:itomoment}, $\EE\|X_t^x\|_H^q<\infty$ for every $q\ge2$, so $P_t^*\delta_x\in\mathcal P_r(H)$ for every $r\ge1$.
\end{proof}

\subsection{Exponential convergence}

It is well known that for $p=2$, with or without additive noise, two solutions driven by the same noise satisfy
\begin{equation}\label{eq:heat}
\|X_t^x-X_t^y\|_H^2\le e^{-2c_0^2t}\|x-y\|_H^2,\quad t>0,
\end{equation}
where $c_0>0$ is the inverse Poincar\'e constant of $U$, that is, $\lambda_1:=c_0^2$ is the first eigenvalue of $-\Delta$ on $U$ with homogeneous Dirichlet boundary conditions. This is a consequence of Poincar\'e's inequality and Gronwall's inequality, and it gives exponential convergence to equilibrium. See~\cite[Chapter 6]{DPZ:96} for the long-time behavior of the stochastic heat equation.
 For small noise intensity and degenerate noise, we obtained in~\cite[Theorem 2.3]{BT:24} the following explicit upper bound for the $\eps$-mixing times, assuming (N1), and that for some $\eta\in (0,1)$,
\[\|B\|^2_{L_2(H_0,H)}\le \eta c_0^2.\]
The bound is stated in $\mathcal W_2$, for unit viscosity, that is, with coefficient $1$ in front of $\Delta$. It reads
\[
\tau_{\mathsf{mix},\mathcal W_2}(\e;x)\leq
\frac{1}{c_0^2}\left[\log\left(\|x\|_H+\frac{\|B\|_{L_2(H_0,H)}}{\sqrt{2(1-\eta)}c_0}\right)+\log\left(\frac{1}{\e}\right)\right].
\]
It is meaningful for $\eps$ so small that its right-hand side is non-negative. Proposition~\ref{prop:heat} below gives a bound of the same form without the smallness and regularity assumptions on $B$.

For $p>2$, in one spatial dimension and for non-degenerate noise, the following exponential convergence holds for the degenerate $p$-Laplace equation. It was proved by Liu~\cite[Theorem 1.5(ii) and Example 3.3]{Liu:09-1} by dimension-free Harnack inequalities for $B=(-\Delta)^{-\theta}$, $\theta\in(\frac14,\frac12]$, and by Wang~\cite[Proposition 6.2]{W:15-1} by asymptotic couplings by reflection under (N3) with $l=1$. There exist constants $C,\lambda>0$ such that
\begin{equation}\label{eq:exp1}\|P_t F-\langle F,\mu_*\rangle\|_\infty\le C\|F\|_\infty e^{-\lambda t},\quad t>0.
\end{equation}
This is an estimate in total variation, uniform over $x\in H$ and over bounded observables, see~\cite[Theorem 1.5(ii)]{Liu:09-1} and~\cite[Theorem 2.2 and Remark 2.2]{W:15-1}. The constants $C,\lambda$ are not explicit. The setting of~\cite{W:15-1} also allows an additional multiplicative noise term, driven by an independent cylindrical Wiener process, whose coefficient is Lipschitz continuous with values in the Hilbert--Schmidt operators. Apart from the well-posedness result of~\cite{GT:16-1}, this is the only result collected here that goes beyond additive noise.

In~\cite[Section 4, Proposition 2]{W:15}, for $d=1$, $p>2$ and purely additive noise satisfying (N3), Wang proves
\[
\|P_tF-\langle F,\mu_*\rangle\|^2_\infty\le C'\left(\langle F^2,\mu_*\rangle-\langle F,\mu_*\rangle^2\right)e^{-\lambda' t},\qquad t\ge1,\quad F\in L^2(\mu_*),
\]
with constants $C',\lambda'>0$ and an explicit lower bound for $\lambda'$ in terms of $p$, $l$ and $\sigma$. From this estimate we derive a rate in $\mathcal W_1$. By Lemma~\ref{lem:moments}, $\mu_*$ has finite second moments. For a Lipschitz function $F$, the variance satisfies, as noted in~\cite{W:15},
\begin{align*}
\langle F^2,\mu_*\rangle-\langle F,\mu_*\rangle^2&=\frac{1}{2}\int_{H\times H}|F(x)-F(y)|^2\,\mu_*(\ud x)\,\mu_*(\ud y)\\
&\le\frac12\Lip(F)^2\int_{H\times H}\|x-y\|_H^2\,\mu_*(\ud x)\,\mu_*(\ud y)\\
&\le\Lip(F)^2\int_H\|y\|_H^2\,\mu_*(\ud y).
\end{align*}
Let $F$ be bounded with $\Lip(F)\le1$. Since the estimate of~\cite{W:15} holds in the supremum norm over $H$, we obtain for every $x\in H$ and $t\ge1$
\[
|P_tF(x)-\langle F,\mu_*\rangle|\le\Big(C'\int_H\|y\|_H^2\,\mu_*(\ud y)\Big)^{1/2}e^{-\lambda't/2}.
\]
An unbounded $1$-Lipschitz function $F$ is the limit of the bounded $1$-Lipschitz functions $(F\wedge n)\vee(-n)$, and the bound extends to $F$ by dominated convergence, since $P_t^*\delta_x$ and $\mu_*$ have finite first moments. Taking the supremum over all $F$ with $\Lip(F)\le1$, Kantorovich--Rubinstein duality yields exponential convergence in $\mathcal W_1$, uniformly in the initial datum,
\begin{equation}\label{eq:expW1}
\sup_{x\in H}\mathcal W_1(P_t^*\delta_x,\mu_*)\le C''e^{-\lambda' t/2},\qquad t\ge1,
\end{equation}
where $C'':=(C'\int_H\|y\|_H^2\,\mu_*(\ud y))^{1/2}$.

For non-degenerate noise, these results give exponential rates, whereas the rate~\eqref{eq:detpol} for the deterministic equation is polynomial. This is an instance of \emph{stabilization by noise}.

\subsection{Summary of known convergence results}\label{subsec:table}

Table~\ref{tab:rates} records the known upper rates together with the distance in which they hold. In the deterministic rows $\mu_*=\delta_0$, so $\mathcal W_r(\delta_{u_t^x},\delta_0)=\|u_t^x\|_H$ for every $r\ge1$, and no moment condition is needed. In the stochastic $\mathcal W_r$ rows, both $P_t^*\delta_x$ and $\mu_*$ belong to $\mathcal P_r(H)$ for every $r\ge1$ by Lemma~\ref{lem:moments}. The stochastic $\mathcal W_r$ rows based on~\eqref{eq:Liu0} and~\eqref{eq:heat} follow from these bounds for two solutions by Lemma~\ref{lem:disi}. For $p=2$ the invariant measure is Gaussian. The rows using $D_\gamma$ require finite $\gamma$th moments, and those based on~\eqref{eq:SSTmu} also require $U$ to be convex and $R_\gamma(x)<\infty$. The total variation row needs no moments. In the last column, ``uniform'' means that the constant does not depend on $x\in H$, and $C(x)$ indicates the dependence of the constant on the initial datum.

\begin{table}[!t]
\caption{Known convergence rates to equilibrium, $D(P_t^*\delta_x,\mu_*)\le C(x)\rho(t)$. The moment assumptions are stated in Subsection~\ref{subsec:table}.}
\label{tab:rates}
\footnotesize
\setlength{\tabcolsep}{3pt}
\begin{tabular}{@{}>{\raggedright\arraybackslash}p{2.8cm}lll>{\raggedright\arraybackslash}p{5.6cm}@{}}
\hline\noalign{\smallskip}
$p$ & Noise & Distance $D$ & Rate $\rho(t)$ & Validity, initial data, constant, source\\
\noalign{\smallskip}\hline\noalign{\smallskip}
$p>2$ & $0$ & $\mathcal W_r$, $r\ge1$ & $t^{-1/(p-2)}$ & $t>0$, $x\in H$, uniform,~\eqref{eq:detpol}\\
$p=2$ & $0$ & $\mathcal W_r$, $r\ge1$ & $e^{-\lambda_1t}$ & $t\ge0$, $x\in H$, $C(x)=\|x\|_H$,~\eqref{eq:heat}\\
$1<p<2$, $p\ge\frac{2d}{d+2}$ & $0$ & $\mathcal W_r$, $r\ge1$ & $e^{-\lambda(x)t}$ & $t>0$, $x\in V_p$, $C(x)=\|x\|_H$, $\lambda(x)=\frac C2\|\nabla x\|_{L^p}^{p-2}$,~\eqref{eq:detexp}\\
$1<p<\frac{2d}{d+2}$, $p\ge2-\frac4d$ & $0$ & $D_\gamma$, $\gamma\le\gamma_*$ & $t^{-p/(2-p)}$ & $t\ge1$, $U$ convex, $x\in W^{1,2}_0$, $C(x)=C_1R_\gamma(x)$,~\eqref{eq:SSTmu}\\
\noalign{\smallskip}
$p>2$ & (N0) & $\mathcal W_r$, $r\ge1$ & $t^{-1/(p-2)}$ & $t>0$, $x\in H$, uniform,~\eqref{eq:Liu0}\\
$p=2$ & (N0) & $\mathcal W_r$, $r\ge1$ & $e^{-\lambda_1t}$ & $t\ge0$, $x\in H$, $C(x)=M_r(x)$,~\eqref{eq:heat}, Prop.~\ref{prop:heat}\\
$\sqrt2\le p<2$, $p>\frac{2d}{d+2}$ & (N0) & $\mathcal W_1$ & $t^{-1/2}$ & $t\ge1$, $x\in H$, $C(x)\lesssim(1+\|x\|_H)^{2/p}$,~\mbox{\cite[Thm.~1.4]{LT:11}}\\
$1<p<\sqrt2$, $p>\frac{2d}{d+2}$ & (N0) & $D_{p^2/2}$ & $t^{-p^2/4}$ & $t\ge1$, $x\in H$, $C(x)\lesssim(1+\|x\|_H)^{p}$,~\mbox{\cite[Thm.~1.4]{LT:11}}\\
$1<p<2$, $p\ge2-\frac4d$ & (N1) & $D_\gamma$, $\gamma\le\gamma_*$ & $t^{-p/(2-p)}$ & $t\ge1$, $U$ convex, $x\in W^{1,2}_0$, $C(x)=C_1R_\gamma(x)$,~\eqref{eq:SSTmu}\\
\noalign{\smallskip}
$p>2$ & (N3) & $D_{\rm TV}$ & $e^{-\lambda t}$ & $d=1$, $l=1$, $t>0$, $x\in H$, uniform,~\eqref{eq:exp1}\\
$p>2$ & (N3) & $\mathcal W_1$ & $e^{-\lambda' t/2}$ & $d=1$, $t\ge1$, $x\in H$, uniform,~\eqref{eq:expW1}\\
\noalign{\smallskip}\hline\noalign{\smallskip}
\end{tabular}
\end{table}

Here $M_r(x):=(\int_H\|x-y\|_H^r\,\mu_*(\ud y))^{1/r}$, and $C(x)\lesssim(1+\|x\|_H)^\alpha$ means $C(x)\le C(1+\|x\|_H)^\alpha$ with $C$ depending on $p,d,U,B$. These exponents follow from the estimates for observables of~\cite{LT:11} stated above for $t\ge1$. Rates in different distances are not comparable without further arguments. For instance, the $D_{p^2/2}$ row gives no information in $\mathcal W_1$. In the deterministic singular rows one may instead use that the solution vanishes after a finite time, by~\cite[Chapter~VII, Section~2]{DB:93} and Proposition~\ref{prop:extinction} for $1<p<2$ with $p\ge\frac{2d}{d+2}$, and by~\cite[Theorems 1.5 and 1.6]{P:11} with $m=\frac{2d}{d-2}$, $d\ge3$, for $1<p<\frac{2d}{d+2}$ with $p\ge2-\frac4d$ and $x\in W^{1,2}_0(U)$, since $W^{1,2}_0(U)\subset L^{2d/(d-2)}(U)$ and $2d/(d+\frac{2d}{d-2})=2-\frac4d$.

\section{Mixing times}\label{sec:mixing}

In this section, we prove the monotonicity of the distance to equilibrium, state the inversion lemma, and derive bounds for the mixing times defined in~\eqref{eq:mixdef}.

\subsection{Monotonicity}\label{subsec:monotone}

\begin{lemma}[Contraction and monotonicity]\label{lem:monotone}
Let $p>1$ and assume (N0). If $1<p<\frac{2d}{d+2}$, so that SVI solutions are used, assume in addition (N1) and that $U$ is convex. Let $X^x$ and $X^y$ be the solutions of~\eqref{eq:pLaplace} with initial data $x,y\in H$, driven by the same Wiener process. Then
\begin{equation}\label{eq:contraction}
\|X_t^x-X_t^y\|_H\le\|x-y\|_H\quad\text{for all }t\ge0,\quad\P\text{-a.s.}
\end{equation}
Consequently, for every $r\ge1$, $s\ge0$ and $\mu,\nu\in\mathcal P_r(H)$,
\begin{equation}\label{eq:nonexp}
\mathcal W_r(P_s^*\mu,P_s^*\nu)\le\mathcal W_r(\mu,\nu),
\end{equation}
and the same holds for $D_\gamma$, $0<\gamma\le1$. If $\mu_*\in\mathcal P_r(H)$ is invariant and $P_t^*\delta_x\in\mathcal P_r(H)$ for all $t\ge0$, then the maps $t\mapsto\mathcal W_r(P_t^*\delta_x,\mu_*)$ and $t\mapsto D_\gamma(P_t^*\delta_x,\mu_*)$ are non-increasing.
\end{lemma}

\begin{proof}
\textit{Step 1. Variational solutions.} First let $p\ge\frac{2d}{d+2}$. Then $X^x$ and $X^y$ are variational solutions in the Gelfand triple $V_p\subset H\subset V_p^*$, see~\cite{PR:07} for $p\ge2$ and~\cite[Section~1]{Liu:09} for $1<p<2$. Write $A(v):=\nabla\cdot(|\nabla v|^{p-2}\nabla v)$ for the drift. Both solutions are driven by the same Wiener process and the noise is additive, so the noise cancels in the difference $w_t:=X_t^x-X_t^y$. Hence
\[
w_t=x-y+\int_0^t\big(A(X_s^x)-A(X_s^y)\big)\,\ud s\qquad\text{in }V_p^*.
\]
The chain rule~\cite[Theorem 4.2.5]{LR:15}, applied to this equation without a noise term, gives
\[
\|w_t\|_H^2=\|x-y\|_H^2-2\int_0^t\int_U\left(|\nabla X_s^x|^{p-2}\nabla X_s^x-|\nabla X_s^y|^{p-2}\nabla X_s^y\right)\cdot\nabla w_s\,\ud z\,\ud s .
\]
The map $a\mapsto|a|^{p-2}a$ is monotone on $\R^d$, that is, $(|a|^{p-2}a-|b|^{p-2}b)\cdot(a-b)\ge0$ for all $a,b\in\R^d$. Hence the double integral is non-negative, and~\eqref{eq:contraction} follows. For $p\ge2$, this argument is contained in the proof of~\cite[Proposition~2.2]{BDP:06}.

\textit{Step 2. SVI solutions.} Now let $1<p<\frac{2d}{d+2}$, and assume (N1) and convexity of $U$. Then the solution is the generalized solution of~\cite[Theorem 4.6]{GT:14}, see~\cite[Definition 2.4, Theorem 2.5 and Corollary 2.7]{SST:23}. Convexity enters only through this existence result, see~\cite[Lemma 2.6]{SST:23}, and is not used in the argument below. It is an SVI solution by~\cite[Appendix C]{GT:14}. By definition, it is the limit in $L^2(\Omega;C([0,T];H))$ of solutions of approximating equations in the Gelfand triple $W^{1,2}_0(U)\subset H\subset W^{-1,2}(U)$. These approximations have initial data $x_m\in W^{1,2}_0(U)$ with $x_m\to x$ in $H$, and noise $B_m\in L_2(H_0,W^{1,2}_0(U))$ with $B_m\to B$ in $L_2(H_0,H)$. We use the same $B_m$ for both initial data $x$ and $y$. Since $|\nabla v|^{p-1}\le1+|\nabla v|$ for $p<2$, the argument of Step~1 applies in this Gelfand triple. It shows that the approximations with initial data $x_m$ and $y_m$ satisfy~\eqref{eq:contraction} with $x_m-y_m$ on the right-hand side. Passing to an almost surely convergent subsequence, we obtain~\eqref{eq:contraction} for $X^x$ and $X^y$, compare~\cite[(3.3)]{GT:16} for Neumann boundary conditions. The limit has continuous paths in $H$, and the Markov property holds by~\cite[Section 6.4]{GT:14}, see~\cite[Subsection 2.2]{SST:23}.

\textit{Step 3. Contraction of transport distances.} Let $\mu,\nu\in\mathcal P_r(H)$, let $\Gamma$ be an optimal coupling of $\mu$ and $\nu$, and let $(\xi,\eta)$ be a random variable with law $\Gamma$, independent of $W$. By the Markov property~\cite[Proposition 4.3.5]{LR:15} in the variational case, and by Step~2 in the SVI case, together with the continuous dependence~\eqref{eq:contraction} on the initial datum, $X_s^\xi$ has law $P_s^*\mu$ and $X_s^\eta$ has law $P_s^*\nu$. Thus $(X_s^\xi,X_s^\eta)$ is a coupling of $P_s^*\mu$ and $P_s^*\nu$, and~\eqref{eq:contraction} gives
\[
\mathcal W_r(P_s^*\mu,P_s^*\nu)^r\le\EE\|X_s^\xi-X_s^\eta\|_H^r\le\EE\|\xi-\eta\|_H^r=\mathcal W_r(\mu,\nu)^r.
\]
The same argument with the cost $\|a-b\|_H^\gamma$ gives the corresponding inequality for $D_\gamma$.

\textit{Step 4. Monotonicity.} Let $s,t\ge0$. By the semigroup property, $P_{t+s}^*\delta_x=P_s^*(P_t^*\delta_x)$, and by invariance, $\mu_*=P_s^*\mu_*$. Hence~\eqref{eq:nonexp} with $\mu=P_t^*\delta_x$ and $\nu=\mu_*$ gives
\[
\mathcal W_r(P_{t+s}^*\delta_x,\mu_*)
=\mathcal W_r(P_s^*P_t^*\delta_x,P_s^*\mu_*)
\le\mathcal W_r(P_t^*\delta_x,\mu_*),
\]
and likewise for $D_\gamma$.
\end{proof}

\begin{remark}\label{rem:markov}
The monotonicity of $t\mapsto\mathcal W_r(P_t^*\delta_x,\mu_*)$ does not follow from the Markov property alone. It rests on the Wasserstein contraction~\eqref{eq:nonexp}, which in turn uses the pathwise contraction~\eqref{eq:contraction}, together with the invariance of $\mu_*$. For total variation, in contrast, the Markov property and invariance suffice. Indeed, every Markov kernel is a contraction in $D_{\rm TV}$, so $D_{\rm TV}(P_{t+s}^*\delta_x,\mu_*)\le D_{\rm TV}(P_t^*\delta_x,\mu_*)$. For multiplicative noise,~\eqref{eq:contraction} fails in general.
\end{remark}

More generally, for a distance or divergence~$\mathsf{dist}$ between probability measures on $H$,\footnote{For instance, the \emph{total variation distance}, the \emph{bounded Lipschitz distance}, the \emph{$\Wcal_r$-Wasserstein distance}, the \emph{Kullback--Leibler divergence}, etc.} we define $\tau_{\mathsf{mix},\mathsf{dist}}(\e;x):=\inf\{t\geq 0: \mathsf{dist}(P^*_t\delta_x,\mu_\ast)\leq \e\}$. This notion is natural when
\begin{itemize}
\item[(i)] $\lim_{t\to \infty}\mathsf{dist}(P^*_t\delta_x,\mu_\ast)=0$, and
\item[(ii)] the map $t\mapsto \mathsf{dist}(P^*_t\delta_x,\mu_\ast)$ is non-increasing as $t$ increases,
\end{itemize}
since then the distance does not exceed $\eps$ after the mixing time. For $\mathsf{dist}=\mathcal W_r$ or $D_\gamma$, property (ii) is provided by Lemma~\ref{lem:monotone} under its hypotheses, and for $D_{\rm TV}$ by Remark~\ref{rem:markov}.

\subsection{From convergence rates to mixing times}\label{subsec:inversion}

The following elementary lemma converts a convergence rate into a bound on the mixing time.

\begin{lemma}[Inversion]\label{lem:inversion}
Let $\mathsf{dist}$ be any distance or divergence between probability measures on $H$, $x\in H$, and suppose that
\[
\mathsf{dist}(P_t^*\delta_x,\mu_*)\le C\rho(t)\qquad\text{for all }t>t_0,
\]
where $C\ge0$, $t_0\ge0$, and $\rho:[t_0,\infty)\to(0,\infty]$ is continuous and strictly decreasing, finite on $(t_0,\infty)$, with $\rho(t)\to0$ as $t\to\infty$. The value $\rho(t_0)=+\infty$ is allowed, as for $\rho(t)=t^{-a}$ with $t_0=0$. Let $\rho^{-1}:(0,\rho(t_0))\to(t_0,\infty)$ be the inverse of the restriction of $\rho$ to $(t_0,\infty)$, and set $\rho^{-1}(s):=t_0$ for $s\in[\rho(t_0),+\infty]$. Then, for every $\eps>0$, with $\eps/0:=+\infty$,
\begin{equation}\label{eq:inversion}
\tau_{\mathsf{mix},\mathsf{dist}}(\eps;x)\le\rho^{-1}\left(\frac{\eps}{C}\right).
\end{equation}
In particular, if $\rho(t)=t^{-a}$ with $a>0$, then $\tau_{\mathsf{mix},\mathsf{dist}}(\eps;x)\le t_0\vee(C/\eps)^{1/a}$, and if $\rho(t)=e^{-\lambda t}$ with $\lambda>0$, then $\tau_{\mathsf{mix},\mathsf{dist}}(\eps;x)\le t_0\vee\lambda^{-1}\log(C/\eps)$, with $\log0:=-\infty$.
\end{lemma}

\begin{proof}
If $C=0$, the distance vanishes for $t>t_0$. Then $\tau_{\mathsf{mix},\mathsf{dist}}(\eps;x)\le t_0$, and $t_0=\rho^{-1}(\eps/C)$ because $\eps/0=+\infty$.

Let now $C>0$. We distinguish two cases. If $\eps/C<\rho(t_0)$, then $t_\eps:=\rho^{-1}(\eps/C)$ satisfies $t_\eps>t_0$ and $C\rho(t_\eps)=\eps$. The hypothesis at $t=t_\eps$ gives $\mathsf{dist}(P_{t_\eps}^*\delta_x,\mu_*)\le\eps$, hence $\tau_{\mathsf{mix},\mathsf{dist}}(\eps;x)\le t_\eps$. If $\eps/C\ge\rho(t_0)$, then $C\rho(t)<C\rho(t_0)\le\eps$ for every $t>t_0$, so the distance does not exceed $\eps$ at any $t>t_0$, and $\tau_{\mathsf{mix},\mathsf{dist}}(\eps;x)\le t_0$. In both cases $\tau_{\mathsf{mix},\mathsf{dist}}(\eps;x)\le\rho^{-1}(\eps/C)$.

For $\rho(t)=t^{-a}$ we have $\rho^{-1}(s)=s^{-1/a}$ for $s<\rho(t_0)$, and for $\rho(t)=e^{-\lambda t}$ we have $\rho^{-1}(s)=\lambda^{-1}\log(1/s)$ for $s<\rho(t_0)$. This gives the two special cases.
\end{proof}

The inversion requires neither monotonicity of $\mathsf{dist}$ in $t$ nor $\eps$ to be small. All upper bounds on mixing times in Subsection~\ref{subsec:ex} and in Table~\ref{tab:mix} rest on estimates $\mathsf{dist}(P_t^*\delta_x,\mu_*)\le C\rho(t)$, $t>t_0$, with $\rho$ non-increasing. They follow either from Lemma~\ref{lem:inversion} or, in Propositions~\ref{prop:degenerate} and~\ref{prop:extinction}, from the same argument applied directly.

\begin{remark}[Upper envelope]\label{rem:envelope}
If the distance to equilibrium is not known to be monotone in $t$, it might exceed $\eps$ again after the mixing time. One may then use
\[
\overline\tau_{\mathsf{mix},\mathsf{dist}}(\eps;x):=\inf\Big\{t\ge0:\sup_{s\ge t}\mathsf{dist}(P_s^*\delta_x,\mu_*)\le\eps\Big\},
\]
the first time after which the distance does not exceed $\eps$. Always $\tau_{\mathsf{mix},\mathsf{dist}}(\eps;x)\le\overline\tau_{\mathsf{mix},\mathsf{dist}}(\eps;x)$, with equality if $t\mapsto\mathsf{dist}(P_t^*\delta_x,\mu_*)$ is non-increasing. The proof of Lemma~\ref{lem:inversion} also applies to $\overline\tau_{\mathsf{mix},\mathsf{dist}}$. Indeed, since $\rho$ is decreasing, the hypothesis of the lemma gives
\[
\sup_{s\ge t}\mathsf{dist}(P_s^*\delta_x,\mu_*)\le C\rho(t)\qquad\text{for all }t>t_0 .
\]
Hence all upper bounds on mixing times in this paper hold for $\overline\tau_{\mathsf{mix},\mathsf{dist}}$ as well, without the monotonicity of Lemma~\ref{lem:monotone}. The lower bound of Proposition~\ref{prop:heatlower} holds for $\overline\tau_{\mathsf{mix},\mathcal W_r}$, since $\tau_{\mathsf{mix},\mathcal W_r}(\eps;x)\le\overline\tau_{\mathsf{mix},\mathcal W_r}(\eps;x)$.
\end{remark}

The constant $C$ and the rate function may depend on $x$, $p$, $B$, $U$ and the distance. For the polynomial rates of Table~\ref{tab:rates}, apart from the maximum with $t_0$, the rate $t^{-1/(p-2)}$ gives the bound $C^{p-2}\eps^{2-p}$ on $\tau_{\mathsf{mix},\mathsf{dist}}(\eps;x)$, the rate $t^{-p/(2-p)}$ gives $C^{(2-p)/p}\eps^{-(2-p)/p}$, the rate $t^{-1/2}$ gives $C^2\eps^{-2}$, and the rate $t^{-p^2/4}$ gives $C^{4/p^2}\eps^{-4/p^2}$. These upper bounds tend to $+\infty$ as $\eps\downarrow0$ if $C>0$, since $\rho>0$. This does not mean that the mixing time itself is unbounded, compare Proposition~\ref{prop:extinction}.

\begin{remark}[Constants]\label{rem:constants}
We track constants explicitly where the argument is self-contained, namely in Propositions~\ref{prop:heat} and~\ref{prop:degenerate}. The constants in Proposition~\ref{prop:extinction} are explicit as well. The bounds that rest on~\cite{Liu:09-1,LT:11,SST:23,JL:09,W:15-1,W:15} contain constants that are not explicit in these sources, except for the explicit lower bound for the exponent $\lambda'$ in~\cite{W:15}. For them we only record on which quantities the constants depend, namely $p$, $d$, $U$, $B$, the order $r$ or $\gamma$, and the initial datum through $\|x\|_H$, $\|\nabla x\|_{L^p}$ or $R_\gamma(x)$. We make no claim of sharper constants than in the cited sources.
\end{remark}

\subsection{Bounds for the mixing times}\label{subsec:ex}

We treat $p=2$, $p>2$ and $1<p<2$ in turn.

\begin{proposition}[The linear case $p=2$]\label{prop:heat}
Let $p=2$ and assume (N0). The case $B=0$ is allowed. Let $\lambda_1=c_0^2$ be the first eigenvalue of $-\Delta$ with homogeneous Dirichlet boundary conditions. Then $\mu_*$ is the centered Gaussian measure with covariance $Q_\infty:=\int_0^\infty e^{s\Delta}BB^*e^{s\Delta}\,\ud s$, and for every $x\in H$, $r\ge1$ and $t\ge0$,
\begin{equation}\label{eq:heattwo}
\mathcal W_r(P_t^*\delta_x,\mu_*)\le M_r(x)\,e^{-\lambda_1t},\qquad M_r(x):=\Big(\int_H\|x-y\|_H^r\,\mu_*(\ud y)\Big)^{1/r}.
\end{equation}
For $r\le2$, $M_r(x)\le M_2(x)=(\|x\|_H^2+\operatorname{tr}Q_\infty)^{1/2}\le(\|x\|_H^2+\|B\|^2_{L_2(H_0,H)}/(2\lambda_1))^{1/2}$. Consequently, for every $\eps>0$,
\begin{equation}\label{eq:heatmix}
\tau_{\mathsf{mix},\mathcal W_r}(\eps;x)\le\frac{1}{\lambda_1}\log_+\left(\frac{M_r(x)}{\eps}\right),\qquad\log_+:=\max\{\log,0\}.
\end{equation}
\end{proposition}

\begin{proof}
The solution is the mild solution
\[
X_t^x=e^{t\Delta}x+\int_0^te^{(t-s)\Delta}B\,\ud W_s,
\]
and its invariant measure is $\mu_*=N(0,Q_\infty)$, see~\cite[Chapter 6]{DPZ:96}. Since $\|e^{s\Delta}\|\le e^{-\lambda_1s}$, the trace of $Q_\infty$ satisfies
\[
\operatorname{tr}Q_\infty=\int_0^\infty\|e^{s\Delta}B\|^2_{L_2(H_0,H)}\,\ud s\le\int_0^\infty e^{-2\lambda_1s}\|B\|^2_{L_2(H_0,H)}\,\ud s=\frac{\|B\|^2_{L_2(H_0,H)}}{2\lambda_1}.
\]
Two solutions driven by the same noise satisfy $X_t^x-X_t^y=e^{t\Delta}(x-y)$, hence $\|X_t^x-X_t^y\|_H\le e^{-\lambda_1t}\|x-y\|_H$. Lemma~\ref{lem:disi} with exponent $r$ gives
\begin{align*}
\mathcal W_r(P_t^*\delta_x,\mu_*)^r&\le\int_H\EE\|X_t^x-X_t^y\|_H^r\,\mu_*(\ud y)\\
&\le e^{-r\lambda_1t}\int_H\|x-y\|_H^r\,\mu_*(\ud y)=e^{-r\lambda_1t}M_r(x)^r,
\end{align*}
which is~\eqref{eq:heattwo}. Since $\mu_*$ is centered, $\int_H\langle x,y\rangle_H\,\mu_*(\ud y)=0$, and therefore
\[
M_2(x)^2=\|x\|_H^2+\int_H\|y\|_H^2\,\mu_*(\ud y)=\|x\|_H^2+\operatorname{tr}Q_\infty .
\]
For $r\le2$, Jensen's inequality gives $M_r(x)\le M_2(x)$. Finally, Lemma~\ref{lem:inversion} with $\rho(t)=e^{-\lambda_1t}$, $t_0=0$ and the constant $M_r(x)$ gives~\eqref{eq:heatmix}.
\end{proof}

For $1\le r\le2$ and $\|x\|_H\le R$, the bound~\eqref{eq:heatmix} gives $\tau_{\mathsf{mix},\mathcal W_r}(\eps;x)\le\lambda_1^{-1}\log_+\big((R^2+\|B\|^2_{L_2(H_0,H)}/(2\lambda_1))^{1/2}/\eps\big)$. Compared with~\cite[Theorem 2.3]{BT:24}, no smallness or spatial regularity of $B$ is needed, and $M_2(x)\le\|x\|_H+\|B\|_{L_2(H_0,H)}/(\sqrt2c_0)$.

The following lower bound for $p=2$ is elementary. It shows that the order $\log(1/\eps)$ in~\eqref{eq:heatmix} cannot be improved.

\begin{proposition}[A lower bound for $p=2$]\label{prop:heatlower}
Let $p=2$ and assume (N0), and let $e_1$ be a normalized eigenfunction of $-\Delta$ with homogeneous Dirichlet boundary conditions for the eigenvalue $\lambda_1$. Then for every $x\in H$, $r\ge1$ and $t\ge0$,
\[
\mathcal W_r(P_t^*\delta_x,\mu_*)\ge|\langle x,e_1\rangle_H|\,e^{-\lambda_1t},
\]
and hence $\tau_{\mathsf{mix},\mathcal W_r}(\eps;x)\ge\lambda_1^{-1}\log(|\langle x,e_1\rangle_H|/\eps)$ for $0<\eps<|\langle x,e_1\rangle_H|$.
\end{proposition}

\begin{proof}
Let $F(y):=\langle y,e_1\rangle_H$. This function is linear and $1$-Lipschitz on $H$. The stochastic integral in the mild solution has mean zero, $e^{t\Delta}$ is self-adjoint and $e^{t\Delta}e_1=e^{-\lambda_1t}e_1$. Hence
\[
\EE F(X_t^x)=\langle e^{t\Delta}x,e_1\rangle_H=\langle x,e^{t\Delta}e_1\rangle_H=e^{-\lambda_1t}\langle x,e_1\rangle_H .
\]
By Proposition~\ref{prop:heat}, $\mu_*$ is centered, so $\langle F,\mu_*\rangle=0$. Kantorovich--Rubinstein duality and $\mathcal W_1(\cdot,\cdot)\le\mathcal W_r(\cdot,\cdot)$ give
\[
\mathcal W_r(P_t^*\delta_x,\mu_*)\ge\mathcal W_1(P_t^*\delta_x,\mu_*)\ge\big|\EE F(X_t^x)-\langle F,\mu_*\rangle\big|=e^{-\lambda_1t}|\langle x,e_1\rangle_H| .
\]
For $t<\lambda_1^{-1}\log(|\langle x,e_1\rangle_H|/\eps)$ the right-hand side is larger than $\eps$, which gives the bound on the mixing time.
\end{proof}

Together with Proposition~\ref{prop:heat}, this gives $\tau_{\mathsf{mix},\mathcal W_r}(\eps;x)=\lambda_1^{-1}\log(1/\eps)+O(1)$ as $\eps\downarrow0$ for every $x$ with $\langle x,e_1\rangle_H\neq0$. For every $a_0>0$, the lower bound of Proposition~\ref{prop:heatlower} is uniform over $\{x\in H:|\langle x,e_1\rangle_H|\ge a_0\}$.

\begin{proposition}[The degenerate case $p>2$]\label{prop:degenerate}
Let $p>2$ and assume (N0). The case $B=0$ is allowed. Let $C_S$ be as in Lemma~\ref{lem:moments}, for instance $C_S=|U|^{\frac12-\frac1p}C_P$, where $C_P$ is the Poincar\'e constant in $\|v\|_{L^p(U)}\le C_P\|\nabla v\|_{L^p(U)}$, $v\in W^{1,p}_0(U)$. Set $\kappa_p:=2^{2-p}C_S^{-p}$.
\begin{enumerate}[(i)]
\item For all $x,y\in H$ and $t>0$, $\P$-a.s., with $0^{2-p}:=+\infty$,
\[
\|X_t^x-X_t^y\|_H\le\left(\|x-y\|_H^{2-p}+(p-2)\kappa_pt\right)^{-\frac{1}{p-2}}\le\left((p-2)\kappa_pt\right)^{-\frac{1}{p-2}}.
\]
\item For all $x\in H$, $r\ge1$ and $t>0$, $\mathcal W_r(P_t^*\delta_x,\mu_*)\le((p-2)\kappa_pt)^{-1/(p-2)}$. Hence, for every $\eps>0$,
\[
\tau_{\mathsf{mix},\mathcal W_r}(\eps;x)\le\frac{\eps^{2-p}}{(p-2)\kappa_p},
\]
uniformly in $x\in H$, $r\ge1$ and $B$.
\end{enumerate}
\end{proposition}

\begin{proof}
(i) Let $w_t:=X_t^x-X_t^y$. As in Step~1 of the proof of Lemma~\ref{lem:monotone}, for a.e.\ $t$,
\[
\frac{\ud}{\ud t}\|w_t\|_H^2=-2\int_U\left(|\nabla X_t^x|^{p-2}\nabla X_t^x-|\nabla X_t^y|^{p-2}\nabla X_t^y\right)\cdot\nabla w_t\,\ud z .
\]
Part~(i) is~\eqref{eq:Liu0} with explicit constants, compare~\cite[Proposition~2.2]{BDP:06} and~\cite[Theorem~1.3]{Liu:09-1}. For $p\ge2$, the inequality $(|a|^{p-2}a-|b|^{p-2}b)\cdot(a-b)\ge2^{2-p}|a-b|^p$ holds for all $a,b\in\R^d$, see~\cite[Section~12, (I)]{Lin:17}. Together with $\|w_t\|_H\le C_S\|\nabla w_t\|_{L^p}$, it gives
\[
\frac{\ud}{\ud t}\|w_t\|_H^2\le-2^{3-p}\|\nabla w_t\|_{L^p}^p\le-2\kappa_p\|w_t\|_H^p .
\]
As long as $w_t\neq0$, this inequality is equivalent to $\frac{\ud}{\ud t}\|w_t\|_H^{-(p-2)}\ge(p-2)\kappa_p$. Integration from $0$ to $t$ gives $\|w_t\|_H^{-(p-2)}\ge\|x-y\|_H^{-(p-2)}+(p-2)\kappa_pt$, which is the first inequality in (i). If $w_{t_1}=0$ for some $t_1$, then $w_t=0$ for all $t\ge t_1$ by~\eqref{eq:contraction}, and (i) holds trivially for these $t$. The second inequality in (i) follows from $\|x-y\|_H^{2-p}\ge0$.

(ii) By Lemma~\ref{lem:moments}, $\mu_*\in\mathcal P_r(H)$. The second bound in (i) does not depend on $y$. Lemma~\ref{lem:disi} with exponent $r$ therefore gives
\[
\mathcal W_r(P_t^*\delta_x,\mu_*)^r\le\int_H\EE\|X_t^x-X_t^y\|_H^r\,\mu_*(\ud y)\le((p-2)\kappa_pt)^{-r/(p-2)} .
\]
The right-hand side is decreasing in $t$, and it is at most $\eps^r$ if and only if $t\ge\eps^{2-p}/((p-2)\kappa_p)$. This gives the bound on the mixing time.
\end{proof}

Under the non-degeneracy assumption (N3) with $l=1$, the mixing times in $D_{\rm TV}$ and, by Remark~\ref{rem:Wr}, in every $\mathcal W_r$ are at most of order $\log(1/\eps)$, uniformly in $x$. Whether degenerate noise can improve the order $\eps^{2-p}$ in $\mathcal W_r$ is open.

\begin{remark}[Non-degenerate noise in one dimension]\label{rem:Wr}
Let $p>2$ and assume (N3) with $l=1$. We show that the total variation bound~\eqref{eq:exp1} implies exponential convergence in every $\mathcal W_r$. Let $\mu,\nu\in\mathcal P_{2r}(H)$, and let $|\mu-\nu|$ denote the total variation measure of $\mu-\nu$. By~\cite[Theorem 6.15]{villani},
\[
\mathcal W_r(\mu,\nu)^r\le2^{r-1}\int_H\|z\|_H^r\,|\mu-\nu|(\ud z).
\]
By the Cauchy--Schwarz inequality, $|\mu-\nu|\le\mu+\nu$ and $|\mu-\nu|(H)=2D_{\rm TV}(\mu,\nu)$, we obtain
\[
\mathcal W_r(\mu,\nu)^r\le2^{r-1/2}\left(\int_H\|z\|_H^{2r}\,(\mu+\nu)(\ud z)\right)^{1/2}D_{\rm TV}(\mu,\nu)^{1/2}.
\]
We apply this with $\mu=P_t^*\delta_x$ and $\nu=\mu_*$. The moments of $P_t^*\delta_x$ are bounded uniformly in $x\in H$ and $t\ge1$. Indeed, let $Y$ have law $\mu_*$ and be independent of $W$. Then $X_t^Y$ has law $\mu_*$, and Proposition~\ref{prop:degenerate}(i) gives $\|X_t^x-X_t^Y\|_H\le((p-2)\kappa_p)^{-1/(p-2)}$ for $t\ge1$. Hence
\[
\big(\EE\|X_t^x\|_H^{2r}\big)^{1/(2r)}\le((p-2)\kappa_p)^{-1/(p-2)}+\Big(\int_H\|y\|_H^{2r}\,\mu_*(\ud y)\Big)^{1/(2r)},
\]
which is finite by Lemma~\ref{lem:moments}. By~\eqref{eq:exp1}, $2D_{\rm TV}(P_t^*\delta_x,\mu_*)\le Ce^{-\lambda t}$ for all $x\in H$. Together, there are constants $C_r>0$ such that
\[
\sup_{x\in H}\mathcal W_r(P_t^*\delta_x,\mu_*)\le C_re^{-\lambda t/(2r)},\qquad t\ge1 .
\]
Lemma~\ref{lem:inversion} with $t_0=1$ gives $\tau_{\mathsf{mix},\mathcal W_r}(\eps;x)\le1\vee 2r\lambda^{-1}\log(C_r/\eps)$ for every $r\ge1$, uniformly in $x\in H$.
\end{remark}

For the deterministic singular equation, the following result is known. We state it with the constants used in Table~\ref{tab:mix}.

\begin{proposition}[The deterministic singular case]\label{prop:extinction}
Let $B=0$, $1<p<2$ and $p\ge\frac{2d}{d+2}$, and let $C_S>0$ be such that $\|v\|_H\le C_S\|\nabla v\|_{L^p(U)}$ for $v\in W^{1,p}_0(U)$. Then for every $x\in H$ the solution of~\eqref{eq:detpLaplace} satisfies $u_t^x=0$ for all $t\ge C_S^p\|x\|_H^{2-p}/(2-p)$. Moreover, for every $r\ge1$ and $\eps>0$,
\[
\tau_{\mathsf{mix},\mathcal W_r}(\eps;x)\le\frac{C_S^p}{2-p}\left(\|x\|_H^{2-p}-\eps^{2-p}\right)_+\le\frac{C_S^p}{2-p}\|x\|_H^{2-p}.
\]
\end{proposition}

\begin{proof}
The vanishing after a finite time is classical, see~\cite[Chapter~VII, Section~2]{DB:93}. The explicit bounds follow from~\cite[Theorem~2.2 and Corollary~3]{BB:20}. Indeed, equation~\eqref{eq:detpLaplace} is the gradient flow in $H$ of the functional $\Phi(v):=\frac1p\int_U|\nabla v|^p\,\ud z$ for $v\in V_p$ and $\Phi(v):=+\infty$ for $v\in H\setminus V_p$~\cite{Br:73}. This functional is convex, lower semicontinuous and absolutely $p$-homogeneous, and $p\Phi(v)\ge C_S^{-p}\|v\|_H^p$ for all $v\in H$ by the choice of $C_S$. Hence the constant denoted $\lambda_1$ in~\cite{BB:20} satisfies $\lambda_1\ge C_S^{-p}$, and~\cite[Corollary~3]{BB:20} gives
\[
\|u_t^x\|_H^{2-p}\le\big(\|x\|_H^{2-p}-(2-p)C_S^{-p}t\big)_+\qquad\text{for all }t\ge0 .
\]
Since $B=0$, we have $\mu_*=\delta_0$ and $\mathcal W_r(\delta_{u_t^x},\delta_0)=\|u_t^x\|_H$ for every $r\ge1$, which gives the bound on the mixing time.
\end{proof}

Thus, for $B=0$, $1<p<2$ and $p\ge\frac{2d}{d+2}$, the mixing times stay bounded as $\eps\downarrow0$, and the exponential bound obtained from~\eqref{eq:detexp} is weaker in this setting. The upper bound is uniform for $x$ in bounded subsets of $H$. For $B\ne0$, the invariant measure is not $\delta_0$, and this argument does not apply.

\begin{proposition}[The singular case with degenerate noise]\label{prop:singular}
Let $1<p<2$ with $p>\frac{2d}{d+2}$, assume (N0), and let $C_S$ be as in Lemma~\ref{lem:moments}. For $x,y\in H$ and $t>0$ let
\[
Z_t(x,y):=\frac1t\int_0^t\left(\|\nabla X_s^x\|_{L^p}+\|\nabla X_s^y\|_{L^p}\right)^p\ud s .
\]
\begin{enumerate}[(i)]
\item For all $x,y\in H$ and $t>0$, $\P$-a.s.,
\[
\|X_t^x-X_t^y\|_H\le\|x-y\|_H\exp\left(-(p-1)C_S^{-2}\,t\,Z_t(x,y)^{-(2-p)/p}\right).
\]
\item For every $r>0$ there is a constant $C_r$, depending only on $r$, $p$, $C_S$ and $\|B\|_{L_2(H_0,H)}$, such that for all $x,y\in H$ and $t>0$,
\[
\EE\|X_t^x-X_t^y\|_H^r\le C_r\|x-y\|_H^r\,t^{-\frac{p}{2-p}}\left(1+t^{-1}\|x\|_H^2+t^{-1}\|y\|_H^2\right).
\]
\item Let $r\ge1$ and $\kappa>0$, and set $\vartheta:=\max\{1,(2-p)r\kappa/p\}$. There is a constant $C=C(p,d,U,B,r,\kappa)$ such that for all $x\in H$ and $t\ge1$,
\begin{equation}\label{eq:singularW}
\mathcal W_r(P_t^*\delta_x,\mu_*)\le C(1+\|x\|_H)^{1+\frac{2\vartheta}{r}}\,t^{-\kappa}.
\end{equation}
Consequently, for all $x\in H$ and $\eps>0$,
\[
\tau_{\mathsf{mix},\mathcal W_r}(\eps;x)\le\max\left\{1,\left(\frac{C(1+\|x\|_H)^{1+2\vartheta/r}}{\eps}\right)^{1/\kappa}\right\}.
\]
In particular, $\tau_{\mathsf{mix},\mathcal W_r}(\eps;x)=O(\eps^{-\delta})$ as $\eps\downarrow0$ for every $\delta>0$.
\end{enumerate}
\end{proposition}

\begin{proof}
(i) Let $w_t:=X_t^x-X_t^y$. As in Step~1 of the proof of Lemma~\ref{lem:monotone}, for a.e.\ $t$,
\[
\frac{\ud}{\ud t}\|w_t\|_H^2=-2\int_U\left(|\nabla X_t^x|^{p-2}\nabla X_t^x-|\nabla X_t^y|^{p-2}\nabla X_t^y\right)\cdot\nabla w_t\,\ud z .
\]
We first bound the integrand from below. By~\cite[Section~12]{Lin:17}, for $\xi,\eta\in\R^d$,
\[
(|\xi|^{p-2}\xi-|\eta|^{p-2}\eta)\cdot(\xi-\eta)\ge(p-1)|\xi-\eta|^2\int_0^1|\eta+s(\xi-\eta)|^{p-2}\,\ud s .
\]
Since $|\eta+s(\xi-\eta)|\le|\xi|+|\eta|$ for $s\in[0,1]$ and $p-2<0$, this gives
\[
(|\xi|^{p-2}\xi-|\eta|^{p-2}\eta)\cdot(\xi-\eta)\ge(p-1)|\xi-\eta|^2(|\xi|+|\eta|)^{p-2}.
\]
We apply this with $\xi=\nabla X_t^x$ and $\eta=\nabla X_t^y$, and write $g_t:=|\nabla X_t^x|+|\nabla X_t^y|$. This gives
\[
\frac{\ud}{\ud t}\|w_t\|_H^2\le-2(p-1)\int_U|\nabla w_t|^2g_t^{p-2}\,\ud z .
\]
Next we compare the right-hand side with $\|\nabla w_t\|_{L^p}$. We write
\[
|\nabla w_t|^p=\big(|\nabla w_t|^2g_t^{p-2}\big)^{p/2}g_t^{(2-p)p/2}
\]
and apply H\"older's inequality with exponents $\frac2p$ and $\frac2{2-p}$. This gives
\[
\|\nabla w_t\|_{L^p}^2\le\|g_t\|_{L^p}^{2-p}\int_U|\nabla w_t|^2g_t^{p-2}\,\ud z .
\]
Let $f_t:=\|\nabla X_t^x\|_{L^p}+\|\nabla X_t^y\|_{L^p}$. Then $\|g_t\|_{L^p}\le f_t$, and since $p-2<0$, also $\|g_t\|_{L^p}^{p-2}\ge f_t^{p-2}$. Together with $\|w_t\|_H\le C_S\|\nabla w_t\|_{L^p}$ we obtain
\[
\frac{\ud}{\ud t}\|w_t\|_H^2\le-2(p-1)f_t^{p-2}\|\nabla w_t\|_{L^p}^2\le-2(p-1)C_S^{-2}f_t^{p-2}\|w_t\|_H^2 .
\]
Gronwall's lemma gives $\|w_t\|_H\le\|x-y\|_H\exp\big(-(p-1)C_S^{-2}\int_0^tf_s^{p-2}\,\ud s\big)$. Finally, the function $s\mapsto s^{(p-2)/p}$ is convex on $(0,\infty)$. Jensen's inequality for the time average therefore gives
\[
\frac1t\int_0^tf_s^{p-2}\,\ud s=\frac1t\int_0^t\big(f_s^p\big)^{(p-2)/p}\,\ud s\ge\Big(\frac1t\int_0^tf_s^p\,\ud s\Big)^{(p-2)/p}=Z_t(x,y)^{-(2-p)/p},
\]
and (i) follows.

(ii) For $s,k>0$ we have $e^{-s}\le(k/e)^ks^{-k}$, since $s\mapsto s^ke^{-s}$ attains its maximum $(k/e)^k$ at $s=k$. We raise (i) to the power $r$ and apply this inequality with $k:=\frac{p}{2-p}$ and $s:=r(p-1)C_S^{-2}tZ_t(x,y)^{-(2-p)/p}$. Since $s^{-k}=(r(p-1)C_S^{-2}t)^{-k}Z_t(x,y)$, taking expectations gives
\[
\EE\|w_t\|_H^r\le\|x-y\|_H^r\left(\frac{k}{e}\right)^k\big(r(p-1)C_S^{-2}t\big)^{-\frac{p}{2-p}}\,\EE Z_t(x,y).
\]
It remains to bound $\EE Z_t(x,y)$. By convexity, $f_s^p\le2^{p-1}(\|\nabla X_s^x\|_{L^p}^p+\|\nabla X_s^y\|_{L^p}^p)$. The energy estimate for $q=2$ in the proof of Lemma~\ref{lem:moments} gives $\EE\int_0^t\|\nabla X_s^x\|_{L^p}^p\,\ud s\le\frac12(\|x\|_H^2+t\|B\|^2_{L_2(H_0,H)})$, and the same for $y$. Hence
\[
\EE Z_t(x,y)\le2^{p-2}\left(2\|B\|^2_{L_2(H_0,H)}+t^{-1}\|x\|_H^2+t^{-1}\|y\|_H^2\right),
\]
which proves (ii).

(iii) We first bound the moments of $Z_t(x,y)$ of every order. Let
\[
\Ecal_t^x:=\frac1t\int_0^t\|\nabla X_s^x\|_{L^p}^p\,\ud s
\]
denote the time-averaged energy. It\^o's formula for $\|X_t^x\|_H^2$, see the proof of Lemma~\ref{lem:moments}, gives
\[
2t\Ecal_t^x\le\|x\|_H^2+t\|B\|^2_{L_2(H_0,H)}+N_t,\qquad N_t=2\int_0^t\langle X_s^x,B\,\ud W_s\rangle_H .
\]
Let $n\ge2$. The quadratic variation of $N$ satisfies $\langle N\rangle_t=4\int_0^t\|B^*X_s^x\|_{H_0}^2\,\ud s\le4\|B\|^2_{L_2(H_0,H)}\int_0^t\|X_s^x\|_H^2\,\ud s$. By the Burkholder--Davis--Gundy inequality, H\"older's inequality in time and~\eqref{eq:itomoment},
\[
\EE|N_t|^n\le C_n\|B\|^n_{L_2(H_0,H)}\,t^{\frac n2-1}\int_0^t\EE\|X_s^x\|_H^n\,\ud s\le C_n\|B\|^n_{L_2(H_0,H)}\,t^{\frac n2}\left(\|x\|_H^n+K'_nt\right).
\]
Dividing the bound on $2t\Ecal_t^x$ by $2t$ and using this estimate, we get, for $t\ge1$,
\[
\EE(\Ecal_t^x)^n\le C\big(\|x\|_H^{2n}+\|B\|^{2n}_{L_2(H_0,H)}+t^{-n/2}(\|x\|_H^n+K'_nt)\big)\le C(1+\|x\|_H^{2n}),
\]
with $C$ depending on $n$, $p$, $C_S$ and $\|B\|_{L_2(H_0,H)}$. For $1\le n<2$ the same bound follows from the case $n=2$ by Jensen's inequality. Since $Z_t(x,y)\le2^{p-1}(\Ecal_t^x+\Ecal_t^y)$, we obtain
\[
\EE Z_t(x,y)^n\le C(1+\|x\|_H^{2n}+\|y\|_H^{2n})\qquad\text{for }t\ge1\text{ and }n\ge1 .
\]
Now we raise (i) to the power $r$ and apply $e^{-s}\le(k/e)^ks^{-k}$ with $k:=r\kappa$ and $s:=r(p-1)C_S^{-2}tZ_t(x,y)^{-(2-p)/p}$. This gives, for $t\ge1$,
\[
\EE\|w_t\|_H^r\le C\|x-y\|_H^r\,t^{-r\kappa}\,\EE\left[Z_t(x,y)^{(2-p)r\kappa/p}\right].
\]
If $(2-p)r\kappa/p\ge1$, the moment bound above applies with $n=(2-p)r\kappa/p$. If $(2-p)r\kappa/p<1$, Jensen's inequality reduces the expectation to the case $n=1$. In both cases
\[
\EE\left[Z_t(x,y)^{(2-p)r\kappa/p}\right]\le C(1+\|x\|_H+\|y\|_H)^{2\vartheta}.
\]
By Lemma~\ref{lem:disi}, by $\|x-y\|_H\le(1+\|x\|_H)(1+\|y\|_H)$ and $1+\|x\|_H+\|y\|_H\le(1+\|x\|_H)(1+\|y\|_H)$, and by the moments of $\mu_*$ of all orders from Lemma~\ref{lem:moments},
\begin{align*}
\mathcal W_r(P_t^*\delta_x,\mu_*)^r&\le\int_H\EE\|X_t^x-X_t^y\|_H^r\,\mu_*(\ud y)\\
&\le Ct^{-r\kappa}(1+\|x\|_H)^{r+2\vartheta}\int_H(1+\|y\|_H)^{r+2\vartheta}\,\mu_*(\ud y).
\end{align*}
The last integral is finite, and taking the $r$th root gives~\eqref{eq:singularW}. The bound on the mixing time follows from Lemma~\ref{lem:inversion} with $\rho(t)=t^{-\kappa}$ and $t_0=1$. For $\delta>0$, the choice $\kappa=1/\delta$ gives $\tau_{\mathsf{mix},\mathcal W_r}(\eps;x)=O(\eps^{-\delta})$.
\end{proof}

Part (ii) contains~\eqref{eq:LT} as the case $r=\frac{2p}{2-p}$. Part (iii) with $r=\frac{2p}{2-p}$ and $\kappa=\frac12$ gives $\vartheta=1$ and the rate $t^{-1/2}$, which is what~\eqref{eq:LT} alone yields. With Proposition~\ref{prop:singular} we generalize the bounds obtained by inverting the estimates for observables of~\cite{LT:11}, namely $\tau_{\mathsf{mix},\mathcal W_1}(\eps;x)\le1\vee(C(x)/\eps)^{2}$ for $p\ge\sqrt2$ and $\tau_{\mathsf{mix},D_{p^2/2}}(\eps;x)\le1\vee(C(x)/\eps)^{4/p^2}$ for $p<\sqrt2$, with $C(x)$ as in Table~\ref{tab:rates}.

Since $D_\gamma(\cdot,\cdot)\le\mathcal W_1(\cdot,\cdot)^\gamma$, Proposition~\ref{prop:singular} also gives rates in $D_\gamma$ faster than any power of $t$ for $1<p<2$, $p>\frac{2d}{d+2}$, under (N0), without convexity of $U$, without (N1) and without the condition $R_\gamma(x)<\infty$. The estimate~\eqref{eq:SSTmu} therefore matters only for $1<p\le\frac{2d}{d+2}$ with $p\ge2-\frac4d$. The disintegration inequality is used with the $r$th power inside the integral. The version without power fails in general for $r>1$, see the discussion after Lemma~\ref{lem:disi}.

We do not know whether the convergence is exponential. By part (i), exponential convergence would follow from a bound $\P(Z_t(x,y)>\zeta)\le C(1+\|x\|_H+\|y\|_H)^ke^{-ct}$ with constants $\zeta,c,k>0$. Compare Proposition~\ref{prop:extinction} for $B=0$.

\subsection{Summary of mixing-time bounds}\label{subsec:tablemix}

Table~\ref{tab:mix} summarizes the upper bounds on mixing times, with the notation of Table~\ref{tab:rates} and Subsection~\ref{subsec:ex}. In the last column, a proposition or remark number refers to the bound stated there, and an equation number to a bound obtained from that estimate in Table~\ref{tab:rates} by Lemma~\ref{lem:inversion}. Each bound holds for every $\eps>0$. The maximum with $1$ in some bounds reflects the threshold $t_0=1$ in Lemma~\ref{lem:inversion}.

\begin{table}[!t]
\caption{Upper bounds on $\tau_{\mathsf{mix},D}(\eps;x)$, valid for all $\eps>0$. The moment assumptions are stated in Subsection~\ref{subsec:table}.}
\label{tab:mix}
\footnotesize
\setlength{\tabcolsep}{3pt}
\begin{tabular}{@{}>{\raggedright\arraybackslash}p{2.8cm}lllll@{}}
\hline\noalign{\smallskip}
$p$ & Noise & $D$ & Upper bound & Data & Source\\
\noalign{\smallskip}\hline\noalign{\smallskip}
$p>2$ & $0$, (N0) & $\mathcal W_r$ & $\eps^{2-p}/((p-2)\kappa_p)$ & $x\in H$ & Prop.~\ref{prop:degenerate}\\
$p=2$ & $0$, (N0) & $\mathcal W_r$ & $\lambda_1^{-1}\log_+(M_r(x)/\eps)$ & $x\in H$ & Prop.~\ref{prop:heat}\\
$1<p<2$, $p\ge\frac{2d}{d+2}$ & $0$ & $\mathcal W_r$ & $C_S^p(\|x\|_H^{2-p}-\eps^{2-p})_+/(2-p)$ & $x\in H$ & Prop.~\ref{prop:extinction}\\
$1<p<\frac{2d}{d+2}$, $p\ge2-\frac4d$ & $0$ & $D_\gamma$ & $1\vee(C_1R_\gamma(x)/\eps)^{(2-p)/p}$ & $x\in W^{1,2}_0$ &\eqref{eq:SSTmu}\\
\noalign{\smallskip}
$1<p<2$, $p>\frac{2d}{d+2}$ & (N0) & $\mathcal W_r$ & $1\vee\big(C(1+\|x\|_H)^{1+2\vartheta/r}/\eps\big)^{1/\kappa}$ & $x\in H$ & Prop.~\ref{prop:singular}\\
$1<p<2$, $p\ge2-\frac4d$ & (N1) & $D_\gamma$ & $1\vee(C_1R_\gamma(x)/\eps)^{(2-p)/p}$ & $x\in W^{1,2}_0$ &\eqref{eq:SSTmu}\\
\noalign{\smallskip}
$p>2$ & (N3) & $D_{\rm TV}$ & $\lambda^{-1}\log_+(C/\eps)$ & $x\in H$ &\eqref{eq:exp1}\\
$p>2$ & (N3) & $\mathcal W_1$ & $1\vee2\lambda'^{-1}\log(C''/\eps)$ & $x\in H$ &\eqref{eq:expW1}\\
$p>2$ & (N3) & $\mathcal W_r$ & $1\vee2r\lambda^{-1}\log(C_r/\eps)$ & $x\in H$ & Rem.~\ref{rem:Wr}\\
\noalign{\smallskip}\hline\noalign{\smallskip}
\end{tabular}
\end{table}

Here $r\ge1$ and $0<\gamma\le\gamma_*$. In the row of Proposition~\ref{prop:singular}, the exponent $\kappa>0$ is arbitrary and $\vartheta=\max\{1,(2-p)r\kappa/p\}$. The constant $\kappa_p$ is defined in Proposition~\ref{prop:degenerate}, $M_r(x)$ in Proposition~\ref{prop:heat}, $C_S$ in Lemma~\ref{lem:moments} and Proposition~\ref{prop:extinction}, $C=C(p,d,U,B,r,\kappa)$ in Proposition~\ref{prop:singular}, $C$ and $\lambda$ in~\eqref{eq:exp1}, $C''$ and $\lambda'$ in~\eqref{eq:expW1}, and $C_r$ in Remark~\ref{rem:Wr}. The upper bounds in the rows $p>2$ are uniform in $x\in H$, and the rows based on~\eqref{eq:exp1} require $l=1$. The $D_\gamma$ rows require $U$ to be convex and $R_\gamma(x)<\infty$. For $p>\frac{2d}{d+2}$, we generalize them in Proposition~\ref{prop:singular}, since $D_\gamma(\cdot,\cdot)\le\mathcal W_1(\cdot,\cdot)^\gamma$. In the deterministic $D_\gamma$ row, the vanishing of the solution after a finite time even gives a bound uniform in $\eps$ for $x\in W^{1,2}_0(U)$, see the paragraph after Table~\ref{tab:rates}.

\section{Conclusion and open problems}\label{sec:conclusion}

We have collected the known quantitative convergence results for the stochastic $p$-Laplace equation in Table~\ref{tab:rates} and the resulting mixing-time bounds in Table~\ref{tab:mix}. As explained in Subsection~\ref{subsec:contrib}, the rows of Table~\ref{tab:mix} cited by an equation number are direct consequences of known estimates via Lemma~\ref{lem:inversion}. The rows cited by a proposition or remark number are stated in Subsection~\ref{subsec:ex}. Except for Proposition~\ref{prop:extinction}, which recalls a known result for the deterministic equation, they are, together with Lemma~\ref{lem:monotone} and Proposition~\ref{prop:heatlower}, the new statements listed in Subsection~\ref{subsec:contrib}. We list some concrete open problems.
\begin{enumerate}[(1)]
\item \emph{Mixing times in the singular regime with degenerate noise.} For $1<p<2$, $p>\frac{2d}{d+2}$, under (N0), Proposition~\ref{prop:singular} gives convergence in $\mathcal W_r$ faster than any power of $t$. Is it exponential, as for $p=2$? By Proposition~\ref{prop:singular}(i), this would follow from exponential tail bounds for the time-averaged energy $Z_t$. For $B=0$ the mixing times are bounded by Proposition~\ref{prop:extinction}, so the question concerns the interaction of noise and singular diffusion.
\item \emph{Lower bounds for $p\neq2$.} For $p\neq2$ and $B\neq0$, find lower bounds on $\tau_{\mathsf{mix},\mathcal W_r}(\eps;x)$ that show to what extent the upper bounds of Propositions~\ref{prop:degenerate} and~\ref{prop:singular} are sharp. By Kantorovich--Rubinstein duality and $\mathcal W_1(\cdot,\cdot)\le\mathcal W_r(\cdot,\cdot)$, $\mathcal W_r(P_t^*\delta_x,\mu_*)\ge|\EE F(X_t^x)-\int_HF\,\ud\mu_*|$ for every $1$-Lipschitz function $F\colon H\to\R$. It therefore suffices to find a non-constant $1$-Lipschitz function $F$ for which the convergence of $\EE F(X_t^x)$ to $\int_HF\,\ud\mu_*$ can be bounded from below, as for $F(y)=\langle y,e_1\rangle_H$ and $p=2$ in Proposition~\ref{prop:heatlower}.
\item \emph{Degenerate noise for $p>2$.} Can degenerate noise with $B\neq0$ improve the order $\eps^{2-p}$ of Proposition~\ref{prop:degenerate}(ii) in $\mathcal W_r$?
\item \emph{Wasserstein rates for non-degenerate noise.} Extend~\eqref{eq:expW1} and Remark~\ref{rem:Wr} to $d\ge2$, and obtain rates in $\mathcal W_r$ that do not deteriorate with $r$, with explicit dependence on the non-degeneracy of $B$.
\item \emph{Dimension and regularity restrictions.} Remove the restrictions $p>\frac{2d}{d+2}$ and $p\ge2-\frac4d$, and the spatial regularity assumption (N1), while keeping control of quantities such as $R_\gamma(x)$.
\item \emph{Cutoff and multiplicative noise.} Study the cutoff phenomenon, that is, an abrupt transition to equilibrium within a time window of smaller order than the mixing time, for small noise intensity, in the spirit of~\cite{BHP:23,BHP:24}. Extend the above to multiplicative noise, where the synchronous coupling is no longer pathwise contractive.
\end{enumerate}

\section*{Declarations}
\noindent
\textbf{Acknowledgments.} 
JMT would like to thank the organizers of the wonderful ``Workshop on Irregular Stochastic Analysis'' in Cortona, Italy, in June 2025.
GB would like to express his gratitude to the Center for Mathematical Analysis, Geometry and Dynamical Systems (CAMGSD) and Instituto Superior T\'ecnico (IST), Lisbon, Portugal, 
for all the facilities used during the preparation of this work.

\noindent
\textbf{Funding.}
The research of both authors was partially supported by the European Union's Horizon Europe research and innovation programme under the Marie Sk\l{}odowska-Curie Actions Staff Exchanges (Grant agreement no.~101183168 -- LiBERA, Call: HORIZON-MSCA-2023-SE-01). Also, the research of GB is funded by FCT/Portugal and the Recovery and Resilience Plan (PRR) through projects UID/04459/2025 (DOI: 10.54499/UID/04459/2025) and UID/PRR/04459/2025 (DOI: 10.54499/UID/PRR/04459/2025).

\noindent
\textbf{Disclaimer.}
Funded by the European Union. Views and opinions expressed are however those of the authors only and do not necessarily reflect those of
the European Union or the European Education and Culture Executive Agency (EACEA). Neither the European Union nor EACEA can be held responsible for them.

\noindent
\textbf{Ethical approval.} Not applicable.

\noindent
\textbf{Competing interests.} The authors declare that they have no conflict of interest.

\noindent
\textbf{Authors' contributions.}
All authors have contributed equally to the paper.

\noindent
\textbf{Availability of data and materials.} Data sharing not applicable to this article as no datasets were generated or analyzed during the current study.

\appendix

\section{The disintegration inequality}

We use the following disintegration result.

\begin{lemma}[Disintegration]\label{lem:disi}
Let $r\ge1$ and let $(P_t)_{t\ge0}$ be a Markovian Feller semigroup on $H$. Assume that $\mu_*\in\mathcal P_r(H)$ is invariant and that $P_s^*\delta_z\in\mathcal P_r(H)$ for all $z\in H$ and $s\ge0$. Then for any $x\in H$ and $0<s\le t$,
\begin{equation}\label{eq:uno}
\mathcal W_r(P_t^*\delta_x,\mu_*)^r
\le\int_H\int_H\mathcal W_r(P_s^*\delta_z,P_s^*\delta_y)^r\,
P_{t-s}^*\delta_x(\ud z)\,\mu_*(\ud y).
\end{equation}
In particular,
\begin{equation}\label{eq:dos}
\mathcal W_r(P_t^*\delta_x,\mu_*)
\le\left(\int_H\mathcal W_r(P_t^*\delta_x,P_t^*\delta_y)^r\,\mu_*(\ud y)\right)^{1/r}.
\end{equation}
\end{lemma}

Lemma~\ref{lem:disi} follows from the convexity of the optimal transport cost~\cite[Theorem~4.8]{villani} with the cost $(a,b)\mapsto\|a-b\|_H^r$. It is applied to the decompositions $P_t^*\delta_x=\int_HP_s^*\delta_z\,P_{t-s}^*\delta_x(\ud z)$ and $\mu_*=\int_HP_s^*\delta_y\,\mu_*(\ud y)$, which hold by the Chapman--Kolmogorov equation and by invariance. The map $z\mapsto P_s^*\delta_z$ is measurable by the Feller property.

For $r=2$ this is~\cite[Lemma 6.1]{BT:24}. For $r=1$,~\eqref{eq:dos} reads
\[
\mathcal W_1(P_t^*\delta_x,\mu_*)\le\int_H\mathcal W_1(P_t^*\delta_x,P_t^*\delta_y)\,\mu_*(\ud y).
\]
For $r>1$, the analogous inequality for $\mathcal W_r$ without the $r$th power is false in general. For $P_t=\mathrm{Id}$ it would read $(\int_H\|x-y\|_H^r\,\mu_*(\ud y))^{1/r}\le\int_H\|x-y\|_H\,\mu_*(\ud y)$. By Jensen's inequality, this fails for every $x\in\operatorname{supp}\mu_*$ unless $\mu_*$ is a Dirac measure.

For~\eqref{eq:pLaplace}, a weaker form of~\eqref{eq:dos}, which suffices in Subsection~\ref{subsec:ex}, can also be seen directly by the stationary synchronous coupling. If $Y$ has law $\mu_*$ and is independent of $W$, then $(X_t^x,X_t^Y)$ is a coupling of $P_t^*\delta_x$ and $\mu_*$, and $\mathcal W_r(P_t^*\delta_x,\mu_*)^r\le\int_H\EE\|X_t^x-X_t^y\|_H^r\,\mu_*(\ud y)$.

\section{Details for the example in Subsection~\ref{subsec:motivation}}

This part gives the details for the motivating example in Subsection~\ref{subsec:motivation}. Throughout, $X_t^x$ denotes the solution of~\eqref{eq:SDE} and $\mu_p$ its invariant measure.

\smallskip\noindent\textit{The bound~\eqref{eq:coupling} for $p\ge2$.} For $p\ge2$, the reflection coupling of~\cite{Eberle} yields constants $C,c>0$ such that
\[
\mathcal W_1(P_t^*\delta_x,P_t^*\delta_y)\le Ce^{-ct}|x-y|\qquad\text{for all }t\ge0\text{ and }x,y\in\mathbb R.
\]
Fix $t\ge0$ and $x,y\in\mathbb R$, and let $(X,Y)$ be an optimal coupling of $P_t^*\delta_x$ and $P_t^*\delta_y$ for $\mathcal W_1$. By the Cauchy--Schwarz inequality,
\[
\mathcal W_2(P_t^*\delta_x,P_t^*\delta_y)^2\le\EE|X-Y|^2\le\big(\EE|X-Y|\big)^{1/2}\big(\EE|X-Y|^3\big)^{1/2}.
\]
The first factor on the right-hand side equals $\mathcal W_1(P_t^*\delta_x,P_t^*\delta_y)^{1/2}$. For the second factor we use that, for $p\ge2$, It\^o's formula gives $\sup_{t\ge0}\EE|X_t^x|^3\le C(1+|x|^3)$. Since $X$ and $Y$ have the laws of $X_t^x$ and $X_t^y$, and since $|a-b|^3\le4(|a|^3+|b|^3)$ for $a,b\in\mathbb R$, we obtain $\EE|X-Y|^3\le C(1+|x|^3+|y|^3)$. Combining these estimates gives
\[
\mathcal W_2(P_t^*\delta_x,P_t^*\delta_y)\le Ce^{-ct/4}|x-y|^{1/4}\big(1+|x|^3+|y|^3\big)^{1/4}\le Ce^{-ct/4}\big(1+|x|+|y|\big).
\]
As $|x|\le1+|x|^p$ for $p\ge1$, we have $1+|x|+|y|\le3(1+|x|^p+|y|^p)$, and this is~\eqref{eq:coupling} with $\lambda_p=c/4$.

\smallskip\noindent\textit{The bound on the mixing time.} Lemma~\ref{lem:disi} holds for $r=2$ on $\mathbb R$ by the same argument. Together with~\eqref{eq:coupling} it gives
\[
\mathcal W_2(P_t^*\delta_x,\mu_p)\le\Big(\int_{\mathbb R}\mathcal W_2(P_t^*\delta_x,P_t^*\delta_y)^2\,\mu_p(\ud y)\Big)^{1/2}\le C_pe^{-\lambda_pt}\bar K_p(x),
\]
where $\bar K_p(x)<\infty$ because $\mu_p$ has moments of every order. Lemma~\ref{lem:inversion}, applied with $\rho(t)=e^{-\lambda_pt}$, $t_0=0$ and the constant $C_p\bar K_p(x)$, gives the bound on $\tau_{\mathsf{mix}}(\eps;x)$ stated in Subsection~\ref{subsec:motivation} for $0<\eps<C_p\bar K_p(x)$.

\smallskip\noindent\textit{The case $1<p<2$.} In this case the drift is not strongly dissipative at infinity. The spectral gap still yields exponential convergence of $P_tf$ in $L^2(\mu_p)$, see~\cite[Chapter~4]{BGL:14}. We show that a bound $\mathcal W_2(P_t^*\delta_x,\mu_p)\le Ce^{-\lambda t}K(x)$ for all $t\ge0$ and $x\in\mathbb R$, with constants $C,\lambda>0$, forces
\[
K(x)\ge\frac{x}{2C}\exp\Big(\frac{\lambda}{4}x^{2-p}\Big)\qquad\text{for all }x\ge1.
\]
In particular, $K$ cannot be polynomial. By the argument of the preceding paragraph, the same holds for the factor $K_p$ in a bound of the form~\eqref{eq:coupling}.

The measure $\mu_p$ is symmetric, so its mean is zero. Testing with the $1$-Lipschitz function $z\mapsto z$ and using $\mathcal W_1(\cdot,\cdot)\le\mathcal W_2(\cdot,\cdot)$, we get
\[
\mathcal W_2(P_t^*\delta_x,\mu_p)\ge\mathcal W_1(P_t^*\delta_x,\mu_p)\ge\EE[X_t^x].
\]
We now bound $\EE[X_t^x]$ from below for $x\ge1$. By It\^o's formula, $\frac{\ud}{\ud t}\EE|X_t^x|^2=1-2\EE|X_t^x|^p\le1$, hence $\EE|X_t^x|^2\le x^2+t$. Again by It\^o's formula, and by Jensen's inequality for the concave function $s\mapsto s^{(p-1)/2}$,
\[
\frac{\ud}{\ud t}\EE[X_t^x]=-\EE\big[(X_t^x)^{[p-1]}\big]\ge-\big(\EE|X_t^x|^2\big)^{(p-1)/2}\ge-(x^2+t)^{(p-1)/2}.
\]
Integrating from $0$ to $t$ gives $\EE[X_t^x]\ge x-t(x^2+t)^{(p-1)/2}$. Now take $t=x^{2-p}/4$. Then $t\le x^2$, so $(x^2+t)^{(p-1)/2}\le(2x^2)^{(p-1)/2}\le\sqrt2\,x^{p-1}$, and therefore
\[
t(x^2+t)^{(p-1)/2}\le\frac{\sqrt2}{4}\,x\le\frac x2.
\]
Hence $\EE[X_t^x]\ge x/2$ for $t=x^{2-p}/4$. Inserting this $t$ into the bound $\mathcal W_2(P_t^*\delta_x,\mu_p)\le Ce^{-\lambda t}K(x)$ gives
\[
\frac x2\le\mathcal W_2(P_t^*\delta_x,\mu_p)\le Ce^{-\lambda x^{2-p}/4}K(x).
\]
Multiplying both sides by $C^{-1}e^{\lambda x^{2-p}/4}$, we obtain $K(x)\ge\frac{x}{2C}\exp\big(\frac{\lambda}{4}x^{2-p}\big)$ for all $x\ge1$, which is the claimed lower bound.

\end{document}